\magnification=\magstep1
\voffset -.5in
\documentstyle{amsppt}
\NoBlackBoxes

\def\a{{\Bbb A}}
\def\coker{\operatorname{ coker}}
\def\mult{\operatorname{ mult}}
\def\C{{\Bbb C}}

\def\Der{\operatorname{ Der}}

\def\Id{\operatorname{ Id}}

\def\E{{\Cal E}} 
\def\Ext{{\operatorname{ Ext}}}

\def\F{{\Cal F}}
\def\G{{\Cal G}}
\def\har#1{\smash{\mathop{\hbox to .8 cm{\rightarrowfill}}
\limits^{\scriptstyle#1}_{}}}
\def\olog#1#2{{\Omega _{#1}(\log #2)}}
\def\oLog#1#2#3{{\Omega^#1 _{#2}(\log #3)}}
\def\hlog#1#2{{\widehat{\Omega} _{#1}(\log #2)}}
\def\hLog#1#2#3{{{\widehat \Omega} ^{#1}_{#2}(\log #3)}}
\def\Olog#1#2#3{{\widehat{{#1}^*\Omega} _{#2}(\log #3)}}
\def\OLog#1#2#3#4{{\widehat{{#1}^*\Omega ^{#2} _{#3}}(\log #4)}}
\def\H{{\Cal H}}

\def\hom{{\Cal H}om}

\def\Im{\operatorname{{im}}}

\def\J{{\Cal J}}

\def\L{{\Cal L}}

\def\M{{\Cal M}}
\def\N{{\Bbb N}}
\def\R{{\Bbb R}}

\def\bo{\operatorname{ O}}
\def\op{\oplus}
\def\ot{\otimes}

\def\O{{\Cal O}}
\def\pr{\noindent{\sl Proof. }}

\def\P{{\Bbb P}}

\def\Q{{\Bbb Q}}
\def\q{{\Cal Q}}
\def\ra{\mathop{\longrightarrow}}
\def\raa{\rightarrow}
\def\re{\noindent{\sl  Remark. }}

\def\S{\hat S}
\def\Sing{\operatorname{ Sing}}

\def\Supp{\operatorname{ Supp}}
\def\t{\tilde}

\def\orb{\operatorname{ orb}}
\def\top{\operatorname{ top}}

\def\Weil{\operatorname{ Weil}}

\def\Z{{\Bbb Z}}

\def\today{\ifcase\month \or January \or February \or March \or April
\or May \or June \or July \or August \or September \or October
\or November \or December \fi \space \number\day,  \number \year}
\def\:{\colon}

\topmatter
\title 
Logarithmic orbifold Euler numbers of surfaces with applications
\endtitle
\author Adrian Langer
\endauthor

\rightheadtext{Orbifold Euler number}

\address{ Adrian Langer:
1. Mathematics Institute, Warwick University, Coventry CV4 7AL, U.K.}\endaddress
\email{langer\@maths.warwick.ac.uk}\endemail
\address{2. Instytut Matematyki UW, ul.~Banacha 2, 02--097 Warszawa, Poland}
\endaddress
\email{alan\@mimuw.edu.pl}\endemail

\thanks
During the preparation of this paper the author enjoyed, as a Marie Curie research fellow,  
the hospitality of the University of Warwick. The author was also partially supported by a Polish KBN grant
(contract number 2P03A02216) and by a subsidy of Prof.~A.~Bialynicki-Birula
awarded by the Foundation for Polish Science.
\endthanks


\endtopmatter

\document 
\vskip .3 cm

\heading 0. Introduction \endheading

If $X$ is a quasiprojective variety with only isolated quotient singularities
then one can define an orbifold Euler number of $X$ as
$$e_{\orb}(X)=e_{\top}(X)-\sum _{x\in \Sing X}\left(1-{1\over r(x)}\right),$$
where $r(x)$ is  the order of the local fundamental group around $x$. 
It is well known that the orbifold Euler number 
computes the top orbifold Chern class of the double dual of the sheaf of regular $1$-forms.
Similarly, using local uniformizations, one can introduce an orbifold Euler number for 
a log canonical surface pair consisting of a normal surface and a reduced Weil divisor on it, 
or more generally  for a log canonical surface pair consisting of a normal surface with
a fractional $\Q$-divisor. This number measures a second Chern class of a suitably defined
logarithmic orbifold vector bundle (or a $\Q$-vector bundle etc.).

In this paper we introduce an orbifold Euler number for any surface pair consisting of 
a normal surface $X$ and a $\Q$-divisor $D$ on it and we interpret it as a second Chern class of
a reflexive sheaf of rational $1$-forms with at most log poles along the $\Q$-divisor $D$.
The logarithmic ramification formula for finite maps and the multiplicativity of Chern classes of reflexive
sheaves on normal surfaces imply the ``proportionality theorem'' (see Corollary 3.7),
generalizing some earlier results by Holzapfel (see [Hz]).
A computation of the orbifold Euler number when $X$ is smooth and $D$ is reduced leads to
a simple proof of a necessary condition for the logarithmic comparison theorem (see Section 6). 
This recovers an earlier result by Calder\'on-Moreno, Castro-Jim\'enez, Mond and Narv\'aez-Macarro
(see [CCMN]).

However the main reason to introduce orbifold Euler numbers is the following generalization
of the Bogomolov--Miyaoka--Yau inequality:

\proclaim{Theorem 0.1}
Let $X$ be a normal projective surface with a $\Q$-divisor $D=\sum a_i D_i$, $0\le a_i\le 1$.
Assume that the pair $(X,D)$ is log canonical and a multiple of $K_X+D$ is effective.
Then
$$3e_{\orb}(X,D)\ge (K_X+D)^2.$$
Moreover, if equality holds then $K_X+D$ is nef.
\endproclaim

This theorem answers all the open questions posed by Megyesi in [Me2], 6.2.
It is different from previously known results [Mi1], [Mi2], [Kob], [KNS], [Wa3], [Me1], [Me2], [La3]
in the fact that we do not assume that the coefficients of $D$ are of a special form, e.g., $1-{1\over n}$.  
This condition imposed very strong restrictions on the support of $D$
restricting applications to curves with ``nice'' singularities like nodes in smooth surfaces
(see [Ti], Theorems 2.7--2.9, [LM], Theorem 3).

The global orbifold Euler number $e_{\orb}(X,D)$ of the pair $(X,D=\sum a_iD_i)$
is defined using local orbifold Euler numbers by
$$e_{\orb}(X,D)=e_{\top}(X)- \sum a_i e_{\top}(D_i-\Sing (X,D))+\sum_{x\in \Sing (X,D)} (e_{\orb}(x;X,D)-1).$$
In applications of Theorem 0.1 the following  properties of local orbifold Euler numbers play a crucial role:

\item{(0.2.1)} If $(X,x)$ is not a quotient singularity and $(X,D)$ is log canonical at $x$ then $e_{\orb}(x;X,D)=0$.
\item{(0.2.2)} If $(X,x)$ is a quotient singularity, $\pi \: (\C ^2, 0)\to (X,x)$ is a finite map
unramified away from $x$ and $K_{\C^2}+{D'}=\pi ^* (K_X+D)$ then
$$e_{\orb}(x;X,D)={1\over \deg \pi} \cdot e_{\orb}(0;\C^2,D').$$
\item{(0.2.3)} If $(\C ^2, D)$ is lc at $0$ and $\mult _0 D$ denotes the multiplicity of $D$ in $0$
(i.e., a sum of multiplicities of irreducible components $D_i$ counted with appropriate multiplicities) then 
$$e_{\orb}(0; \C ^2,D)\le  \left( 1-{\mult _0 D \over 2}\right)^2.$$

In particular, these properties imply that $e_{\orb}(x; X,D)\le 1$ for any lc pair $(X,D)$.
As a corollary we also get the following result:

\proclaim{Theorem 0.3}
Let $D=\sum a_iD_i$ be an effective $\Q$-divisor on a smooth projective surface $X$, where $D_i$ are distinct 
irreducible curves on $X$. Assume that $(X,D)$ is log canonical and a multiple of $K_X+D$ is effective.
Let $m_P$ be the multiplicity of $D$ at $P$,  $\Sing D$ the set of all singular points of the support of $D$,
$r_{P,i}$ the number of analytic branches of $D_i$ passing through $P$ and $g(D_i)$ the geometric genus of $D_i$.
Set $r_P=\sum _{P\in \Sing D} a_i r_{P,i}$. Then
$$
(K_X+D)^2\le 3\left( c_2(X)+\sum a_i(2g(D_i)-2) +\sum _{P\in \Sing D} \left(r_P-m_P+{m_P^2\over 4}\right) \right) .$$
\endproclaim

This theorem generalizes [Ti], Theorem 2.7 (its naive generalization was conjectured in [Ti], 
Remark after Theorem 2.7) and [LM], Theorem 4.
Let us also note that $(X,D)$ is canonical, and in particular lc, if 
$ m_P\le 1$ at each singular point $P$ of $D$ (see [Kol], Exercise 3.14).
Moreover, if $(X,D)$ is log canonical (lc) then $m_P\le 2$ (by a simple computation on the 
blow up at $x$; see also [Kol], Lemma 8.10). Equality $m_P=2$  can hold, e.g.,
if the support of $D$ has an ordinary singularity at $x$. 

The inequalities obtained in Theorems 0.1 and 0.3 are sharp in the sense that equality holds for infinitely many
curves with arbitrarily high multiplicity of singular points (see Example 11.3.2). These examples are not 
connected to ball quotients as in the usual Bogomolov--Miyaoka--Yau inequality.
However, even in our case these inequalities should reflect the existence of an approximate singular 
K\"ahler--Einstein metric on 
the universal orbifold covering of the pair $(X,D)$, with equality being equivalent to constant bisectional
curvature (or to flatness of the corresponding Higgs bundle). Our method, being algebraic, gives no
insight into this type of problems.

\medskip
By Theorem 0.1 if $K_X^2=3c_2(X)>0$ and a multiple of $K_X$ is effective then $K_X$ is nef and big.
If we apply Theorem 0.3 for any irreducible curve $C$ then we get
$$2K_XC\le 6g(C)-6+\sum _{P\in \Sing C}3(r_P-m _P) ,$$
where $r _P$ is the number of analytic branches of $C$ at $P$ and $m _P$ is the multiplicity of $C$ at $P$.
Since $r_P\le m_P$ the last inequality implies $K_XC\le 3g(C)-3$. 
In particular, $X$ does not contain $(-2)$-curves and hence $K_X$ is ample.
Moreover, $X$ does not contain rational and elliptic curves (even singular).

More generally, application of our results to curves in surfaces of general type leads to effective versions
of Bogomolov's result on boundedness of rational curves in surfaces of general type
with $c_1^2>c_2$. For precise results and the history of the problem we refer to Section 10.
Another application leads to better than previously known bounds on singularities of plane curves
with simple or ordinary singularities, e.g, for arrangements of lines, conics etc.~(see Section 11).  

\medskip

Theorem 0.1 is obtained in the same way as in [La3] except that we need the logarithmic ramification formula
for all log canonical pairs. This formula (see Theorem 4.9) can be thought of as a generalization 
to  all log canonical pairs of Deligne's 
result about closedness of logarithmic forms for normal crossing divisors on smooth varieties. 
It also implies the Bogomolov--Sommese type vanishing theorem for all log canonical
pairs (see Theorem 4.11), which is the main ingredient of the proof of Theorem 0.1. 
This part of the paper works in any dimension and it is closely related to 
the problem of characterization of log canonical pairs
in terms of symmetric powers of logarithmic forms with poles along a $\Q$-divisor (see 4.1--4.7).

A large part of the paper is devoted to the study of local orbifold Euler numbers.
These invariants of local surface pairs cannot be computed from
the graph of the minimal resolution and therefore they are only analytic and not topological invariants.
However, we can compute these invariants, e.g., if the support of a $\Q$-divisor has only ordinary singularities
(see Section 8) or for log canonical pairs  ``with at most fractional boundary'' (see Section 9).

As we already noted the local orbifold Euler numbers are particularly well behaved for log canonical pairs.
If a surface germ $(X,x)$ is Kawamata log terminal then the log canonical threshold 
of a divisor $D$ on $(X,x)$  should correspond to the minimal $\alpha$ for which the local orbifold Euler number of
$(X, \alpha D)$  at $x$ vanishes. 

\medskip
The structure of the paper is as follows. In Section 1 we introduce local relative Chern numbers
and list some of their properties. Section 2 is devoted to virtual sheaves of logarithmic 1-forms 
with poles along $\Q$-divisors. In Section 3 we introduce local and global orbifold Euler numbers
and we interpret them as top Chern numbers of the sheaf of rational 1-forms with at most log poles.
In Section 4 we prove the logarithmic ramification formula for log canonical pairs and the corresponding
Bogomolov--Sommese type vanishing theorem. This leads to Theorems 0.1 and 0.3 in Section 5. In Section 6
we compute orbifold Euler numbers for reduced divisors on smooth surfaces and apply them to
study logarithmic cohomology of open surfaces. Sections 7, 8 and 9 are devoted to properties 
and to computation of local orbifold Euler numbers. These results are applied in Sections 10 and 11
to study plane curves and curves in surfaces of general type.

\medskip 

Disclaimer: our invariants and their construction have no obvious connection with the stringy orbifold Euler numbers
introduced by Batyrev, Totaro, Borisov and Libgober. In fact, their stringy orbifold Euler numbers 
can be computed from the resolution graph and therefore they are topological invariants of the underlying 
pair.

\medskip

\noindent
{\sl Notation.}

A $\Q$-divisor $D=\sum d_i D_i$ with $D_i$ distinct prime divisors is called a {\sl boundary}
if $0\le d_i \le 1$ for all $i$. If $f\: Y\to X$ is a birational map then $f^{-1}D$ stands for 
the strict transform of $D$.

A resolution of singularities $f\: {\t X}\to X$ is called a {\sl log resolution} of the pair $(X,D)$
if the exceptional divisor and the strict transform of $D$ have normal crossings. Let $E_i$
be the exceptional divisors. Write
$$K_{\t X}+f^{-1}D+\sum E_i=f^*(K_X+D)+\sum a_i E_i.$$
A surface $X$ with boundary $D=\sum d_i D_i$  is called {\sl log canonical} ({\sl log terminal}; {\sl klt})
if for any log resolution $f\: {\t X}\to X$ of $(X,D)$
we have $a_i\ge 0$ ($a_i>0$; $a_i>0$ and $d_i<1$, respectively) for all $i$. 

We say that $(X,D)$ has a {\sl fractional boundary} if all the coefficients of $D$ are of the form $1-1/m$,
where $m\in \N\cup \infty$.
A log canonical pair $(X,D)$ has {\sl at most fractional boundary} if there exists a log canonical pair 
$(X,D')$ with fractional boundary and such that $D'\ge D$.

Finally, $\S ^n \F$  denotes the double dual of the $n$-th symmetric power of $\F$.
 
\medskip

\heading 1. Preliminaries \endheading

In this preliminary section we recall and extend a definition of local relative Chern numbers,
and we show their basic properties. First, we need the following frequently used auxiliary definition.

\proclaim{Definition 1.1}
A \rom{square } is a commutative diagram
$$\CD
 {\t Y} @>{g}>> Y \\
@VV{\t \pi}V        @VV{\pi }V \\
 \t X @>{f}>> X\\ 
\endCD
$$
in which
\item{1.} $X$, $\t X$, $Y$ and $\t Y$ are normal surfaces
(not necessarily complete),
\item{2.} all the maps are proper,
\item{3.} $f$ and $g$ are birational,
\item{4.} $\pi$ and $\t \pi$ are finite.
\endproclaim

Let $(X,x)$ be a germ of a normal surface singularity and 
$f\: (\t X, E)\to (X,x)$ any resolution of $(X,x)$.

\proclaim{Definition 1.2} ([Wa2, (3.1)])
For a vector bundle $\F$ on $\t X$ we define the
\rom{ modified Euler characteristic} by
$$\chi (x, \F)=\dim \left( (f_*\F)^{**}/f_*\F \right) +\dim R^1 f_*\F .$$  
\endproclaim

\proclaim{Definition 1.3} 
Let $\F$ be a rank $2$ vector bundle on $(\t X,E)$.

\noindent (1)
The \rom{ first Chern class} $c_1(x,\F)$ is an $f$-exceptional $\Q$-divisor
whose intersection with any $f$-exceptional divisor $F$ is equal to the
degree of $\F |_F$.

\noindent (2) The \rom{ second RR Chern class} of $\F$ is defined by
$$c_2(x,\F)={1\over 4}c_1(x,\F)^2+{3\over 4}\liminf_{n\to \infty}
{{\chi (x,S^{2n}\F(-n\det \F))}\over{n^3}}.$$
\endproclaim

The notation used in this paper differs from that used in [Wa2], [La1], [La2] and [La3]:
we write $c_2(x, \F )$ instead of $c_2'(x, \F )$. The reason is that  in this paper 
we do not use  Wahl's definition of second local Chern class. 

\medskip

If $f\: Y\to X$ is a birational proper morphism from a smooth surface $Y$
to a normal surface $X$  and $\F$ is a vector bundle on $Y$ then we set
$$\chi (f, \F)= \sum _{x\in \{ y\in X\: \dim f^{-1}(y)>0\}} \chi (x, \F)$$
and
$$c_i (f, \F)= \sum _{x\in \{ y\in X\: \dim f^{-1}(y)>0\}} c_i(x, \F )$$
for $i=1$, $2$.

\proclaim{Lemma 1.4}
Let $g\: Z\to Y$ and $f\: Y\to X$ be morphisms of normal surfaces and $\F$ a reflexive sheaf on
$Z$. Then
$$\chi (fg, \F)=\chi (g, \F)+\chi (f, \G),$$
where $\G= (g_*\F)^{**}$.
\endproclaim

\pr
From the Leray spectral sequence $R^pf_*(R^qg_*\F) => R^{p+q}(fg)_*\F$
one can easily get an exact sequence
$$0\to R^1f_* (g_*\F)\to R^1(fg)_*\F \to f_* (R^1g_*\F)\to 0.$$
On the other hand if $\H =\G /g_*\F$ then
$$0\to f_* (g_*\F)\to f_*\G \to f_* \H \to  R^1f_* (g_*\F)\to R^1f_*\G\to 0.$$
Now the lemma follows from the above sequences by an  easy calculation, Q.E.D.

\medskip

In the following we  need a generalization of relative local Chern classes
to any birational morphism.

\proclaim{Definition--Proposition 1.5}
Let ${\hat f}\: {\hat X}\to X$ be a proper birational map of normal surfaces
and let $\F$ be a reflexive sheaf on $\hat X$. Let 
$\alpha \:{\t X}\to {\hat X}$ be a resolution of singularities and set
$f=\alpha {\hat f}$. Take any vector bundle $\G$ on $\t X$ such that
$(\alpha _* \G) ^{**}=\F$. Then the number 
$$c_2({\hat f}, \F)=c_2(f,\G)-c_2(\alpha, \G)$$
does not depend on the choice of $\alpha$ and $\G$ and we call it 
a \rom{relative second Chern class} of $\F$ with respect to $f$.
\endproclaim

The above proposition follows easily from Lemma 1.4.

\proclaim{Theorem 1.6}
If we have a square
$$\CD
 {\t Y} @>{g}>> Y \\
@VV{\t \pi}V        @VV{\pi }V \\
 \t X @>{f}>> X\\ 
\endCD
$$
then 
$$c_2(({\t \pi}^*\F)^{**},g)=\deg \pi \cdot c_2(\F, f)$$
for any rank $2$ reflexive sheaf $\F$ on $\t X$.
\endproclaim

\pr
The proposition follows easily from the fact that the relative second
Chern class of vector bundles on smooth surfaces behaves
multiplicatively under generically finite proper maps 
(see [La3], Theorem 1.5). Q.E.D.

\medskip

Now we can introduce Chern classes of reflexive sheaves.

\proclaim{Definition--Proposition 1.7} (see
[La2, Definition--Proposition 2.8])
Let $\E$ be a rank $2$ reflexive sheaf on a normal proper surface $X$. 
The first Chern class $c_1\E$ of $\E$ is the class of $\det \E$ in $\Weil X$. 
Let $f\:Y\to X$ be any resolution of 
singularities and $\F$ any vector bundle on $Y$ such that 
$(f_*\F)^{**}=\E$. Then the number
$$c_2 \E=c_2{\F}-c_2(f,\F)$$
does not depend on the choice of $f$ and $\F$ and we call it the
\rom{ second RR Chern class} of $\E$.
\endproclaim

For further properties of the above Chern classes we refer the reader
to [Wa2], [La1] and [La2, Section 2].

\proclaim{Theorem 1.8} (cf.~[Wa2, Theorem 1.8])
Let $\E$ be a rank 2 vector bundle on a smooth projective curve $C$
and let $N$ be a degree $d>0$ line bundle. Set $e=\deg \E$ and
$${\overline s}={\overline s} (\E)=\max \left( {e\over 2}, \max \{\deg \L \:
\L\subset \E \}\right) .$$ 
Then
$$\sum _{i>0} h^0(C, S^{2n}\E (-n\det \E -i N))= {n^3\over 6d}(2{\overline s}-e)^2
+\bo (n^2).$$
\endproclaim

\pr
If $\E$ is semistable then $S^{2n}\E(-n\det \E-iN)$ is semistable 
of negative degree. Hence
$h^0(S^{2n}\E(-n\det \E-iN))=0$ for $i>0$ and the theorem follows.

Now assume that $\E$ is unstable and let $L\subset \E$ be a line subbundle
of maximal degree $\deg L=s$. Set  $M=\E/L$. Taking sequences of symmetric 
powers of the sequence
$$0\raa L\ra \E\ra M\raa 0$$  
we get a filtration of $S^{2n}\E(-n\det \E-iN)$ with quotients
$\{(n-j)(M-L)-iN\}$. Since $\deg(L-M)=2s-e>0$ we get
$$\sum _{i>0}h^0(S^{2n}\E(-n\det \E-iN))\le \sum_{j=0}^{n}\sum_{i>0}
h^0(j(L-M)-iN).$$
By Lemma 1.9, [Wa2], the right hand side of this inequality is equal to
$$\sum_{j=0}^n\left\lceil{j(2s-e)-(2g-2)\over d}\right\rceil\left(j(2s-e)-
{d\over2}\left\lceil{j(2s-e)-(2g-2)\over d}\right\rceil\right)+O(j)=$$
$$
=\sum_{j=1}^nj^2{(2s-e)^2\over 2d}+O(n^2)={n^3\over 6d}(2s-e)^2+O(n^2).
$$

\medskip
On the other hand
$$h^0(S^{2n}\E(-n\det \E-iN))\ge h^0(S^\delta \E((2n-\delta)L-n\det \E-iN)).$$
Set $\delta=n-[{id+2g-1\over 2s-e}]$ for $0< i\le {n(2s-e)-(2g-1)
\over d}$. Then
$$h^0(S^\delta \E((2n-\delta)L-n\det \E-iN))=(\delta+1)
\left( {e\delta\over 2}+(2s-e)n-s\delta-id+1-g\right)=$$
$$=
n^2(s-e/2)-idn+i^2{d^2\over
2(2s-e)}+O(n).$$
Summing over $0< i\le {n(2s-e)-(2g-1)\over d}$ we get
${n^3\over 6d}(2s-e)^2+O(n^2)$. Therefore
$$\sum _{i>0}h^0(S^{2n}\E(-n\det \E-iN))\ge {n^3\over 6d}(2s-e)^2+O(n^2),$$
which completes the proof.

\medskip
This theorem is similar in statement and proof to Theorem 1.8, [Wa2].
It also follows quite easily from this theorem if we use Corollary 4.19, [La1]
and Example 2.8 and Theorem 3.11, [Wa2].
However, in our case a direct proof is more transparent than the original
proof of Theorem 1.8, [Wa2].

As a particular case of the above theorem we get the following characterisation of
rank $2$ semistable bundles:

\proclaim{Corollary 1.9} 
$\E$ is semistable if and only if
$$\sum _{i>0} h^0(C, S^{2n}\E (-n\det \E -i N))= \bo (n^2).$$
\endproclaim

\proclaim{Theorem 1.10} ([Wa2], Example 2.8 and Theorem 3.11)
Let $\t X$ be the total space of a line bundle $N^{-1}$ and let
$\pi \: {\t X}\to C$ be a canonical projection. Let $f\: {\t X}\to X$
be a contraction of the zero section of $N^{-1}$. Then
$$c_2(f,\pi^*\E)=-{{\overline s} (e-{\overline s})\over d}.$$
\endproclaim

\medskip

\heading 2. Logarithmic forms with poles along $\Q$-divisors \endheading 

In this section we introduce virtual sheaves of rational  $q$-forms with log poles along $\Q$-divisors.
They are used in the next section to define local and global orbifold Euler numbers.
\medskip

Let $(X,D)$ be a normal projective variety with a 
$\Q$-divisor $D=\sum a_i D_i$,
where $D_i$ are distinct integral divisors and $0\le a_i\le 1$.
There exists a normal variety $Y$ and a finite dominant morphism
$\pi \: Y\to X$  such that $\pi^*(a_i D_i)$ are integral Weil divisors.
Let $U= X- \Sing(X,D)- \pi (\Sing Y)$ and let $j$ be an embedding 
$V=\pi ^{-1}(U)\hookrightarrow Y$. First we define a vector bundle on $V$ by
setting 
$${\OLog {(\pi|_U)} q U {D|_U}}=\bigwedge ^q\left({\pi^*dz\over t} \O_V +\pi^*\Omega _U\right),$$
where $z=0$ is a local equation of $\Supp D$ and $t=0$ is a local equation
of $\pi^*D$. 

\proclaim{Definition 2.1}
In this set up we define a \rom{pull back of $\O _X$-module of 
logarithmic $q$-forms} of $X$ along $D$
to $Y$ as an $\O_Y$-module
$${\OLog \pi q X D}=j_*({\widehat{(\pi|_U) ^{*} \Omega ^q} _{U}(\log D|_U)}).$$ 
\endproclaim

\proclaim{Lemma 2.2}
${\OLog \pi q X D}$ is a reflexive sheaf.
\endproclaim

\pr
By definition ${\OLog \pi q X D}$ is of the form $j_*\F$ for some vector bundle
$\F$ and an inclusion $j\:V\to Y$ such that $Y-V$ is a closed subset of $Y$ of
codimension $\ge 2$. This implies that the sheaf ${\OLog \pi q X D}$ is normal, i.e., 
one can extend its sections 
defined outside a codimension $\ge 2$ subset. Since it is also torsion free, it must
be reflexive, Q.E.D.

\medskip
We include this simple lemma since the original proof 
of Saito ([Sai2], Corollary 1.7)
works only for $q=1$ and in the literature there still appear proofs of special cases
of the above lemma. Let us note that in much the same way as above one can define
a sheaf ${\widehat{\pi ^{*} \Der} _{X}(\log D)}$ of {\sl pull backs of logarithmic vectors along
the $\Q$-divisor} $D$. Then ${\OLog \pi q X D}\simeq \hom _{\O _Y}(\bigwedge ^q 
{{\widehat {\pi ^{*} \Der} _{X}(\log D)}}, \O_Y)$, since both sheaves are reflexive
and they are isomorphic in codimension $1$.

Now let us change to the analytic category, which is completely analogous.

\medskip
Let $U$ be a domain in $\C ^n$ and let $D$ be a $\Q$-divisor on $U$
with support defined by an equation $h=0$, where $h$ is holomorphic 
on $U$. Let $\pi \: V\to U$ be a finite holomorphic mapping from
a normal analytic set $V$ such that $\pi^*D$ is a Cartier divisor
given by $t=0$.

\proclaim{Lemma 2.3}(cf. [Sai2], (1.1))
Let $\omega$ be a meromorphic $q$-form on $V$ (with poles only along $\pi^*D$).
Then the following conditions are equivalent:
\item{(1)} $\omega \in {\OLog \pi q U D}$,
\item{(2)} $t\omega \in \pi^* \Omega ^q_U$ and $\omega \wedge \pi^*(dh)\in \pi^*\Omega^{q+1} _U$,
\item{(3)} $t\omega \in \pi^* \Omega _U^q$ and $ \pi^*(h)\cdot d\omega\in \pi^*\Omega _U ^{q+1}$,
\item{(4)} There exist a holomorphic function $g$ on $V$,  $\xi \in \pi^* \Omega ^{q-1}_U$  and 
$\eta \in \pi^* \Omega _U ^q$ such that the set $\{ z\in V \: g(z)=t(z)=0\}$ has
codimension at least $2$ and
$$g\omega ={\pi^*(dh)\over t}\wedge \xi +\eta ,$$
\item{(5)} There exists an $(n-2)$-dimensional analytic set $A\subset \Supp (\pi^*D)$ such that
the germ of $\omega$ at any point $P \in \Supp (\pi ^*D)-A$ belongs to $({\pi^*(dh)\over t}
\wedge\pi^*\Omega _U^{q-1} +\pi^*\Omega _U^q)_P$.
\endproclaim

We skip the proof, which is analogous to the proof of [Sai2], (1.1).

\medskip

\heading 3. Orbifold Euler numbers \endheading

In this section we introduce local and global orbifold Euler numbers. Then we prove that they
compute top Chern numbers of corresponding sheaves 
of logarithmic 1-forms with poles (Theorem 3.6)
and use this to prove the ``proportionality theorem'' (Corollary 3.7).

\medskip

Let $(X,x)$ be a germ of a normal surface singularity and 
let $D=\sum a_i D_i$  be a boundary $\Q$-divisor on $X$. 
Let $f\:(\t X, E)\to (X,x)$ be a log resolution of the pair $(X,D)$ at $x$,
where $E$ denotes a reduced scheme structure on the exceptional set of $f$.
Let us assume that there exists a finite proper map $\pi\: (Y,y)\to (X,x)$ from a normal
surface $Y$ such that 
$\pi^*D$ is a Weil divisor (note that $y=\pi^{-1}(x)$ need not be a point).   
Then we can construct a square
$$\CD
 (\t Y,F) @>{g}>> (Y,y) \\
@VV{\t \pi}V        @VV{\pi }V \\
 (\t X, E) @>{f}>> (X,x).\\ 
\endCD
$$

If $(X,D,x)$ is a germ of a quasi-projective surface with boundary then
such a map $\pi$ always exists and  we can assume that $\t Y$ is smooth (see Lemma 2.5, [La3]). 
Indeed, by Kawamata's covering lemma (see [KMM, Theorem 1-1-1]) there exists
a finite Galois covering ${\t \pi}\: {\t Y}\to {\t X}$ such that
$\t Y$ is smooth and ${\t \pi} ^{-1}(D+E)$ is a normal crossing  Weil divisor.
Then $g\:{\t Y}\to Y$ and $\pi\:Y\to X$ come from the Stein factorization of 
the composition $f {\t \pi}\: {\t Y}\to X$. 

In the analytic case if we have a finite surjective map $\pi$ between normal surface
singularities then taking an appropriate resolution 
$f$ we can assume that $\t Y$
has at most cyclic quotient singularities.

\proclaim {Definition 3.1}
A \rom{local orbifold Euler number of the pair $(X,D)$} at $x$ is defined by
$$e_{\orb} (x;X,D)=-{c_2(g, {\Olog {\t \pi} {\t X} {f^{-1}D+E}})\over \deg
\pi }.$$
\endproclaim

\noindent
{\sl 3.2.} We will show that $e_{\orb }(x;X,D)$ is well defined and it depends only 
on the analytic type of $(X,D)$ at $x$. 

First note that this number does not depend on the choice of  
$\pi$. Indeed, if we have two maps 
$\pi_1\:(Y_1,y_1)\to (X,x)$ and $\pi_2\: (Y_2,y_2)\to (X,x)$ then
we can construct $\pi_3\:(Y_3,y_3)\to (X,x)$ by taking the normalization
of the fiber product $Y_1\times _{X} Y_2$. 
Set $\F_i={\Olog {{\t \pi} _i} {\t X} {f^{-1}D+E}}$ for $i=1,2,3$.
Then 
$$c_2 (g_2; \F_2)=
{1\over \deg \pi _1}c_2 (g_3; \F_3)={\deg\pi_2\over \deg \pi _1}
c_2 (g_1; \F_1),
$$
which proves our claim.

Now take any two log resolutions of the pair $(X,D)$ at the point $x$.
Since any two log resolutions of $(X,D)$ are dominated 
by a third one we can assume that one resolution
$f_2\:({\t X}_2, E_2)\to (X,x)$ dominates the other resolution
$f_1\:({\t X}_1, E_1)\to (X,x)$. Actually, it is sufficient to prove
the required equality of Chern classes if the map
$\alpha\:({\t X}_2, E_2)\to ({\t X}_1, E_1)$ is a blow up at a single point. 
Let $\beta\:({\t Y}_2, F_2)\to ({\t Y}_1, F_1)$ be the induced map and 
set $\F_i={\Olog {{\t \pi} _i} {\t X} {f_i^{-1}D+E}}$ for $i=1,2.$
Then $(\beta _*\F_2)^{**}=\F_1.$ Hence to prove that $c_2(g_1,\F_1)=c_2(g_2,\F_2)$
it is sufficient to show that 
$c_2(\beta, \F_2)=0$. This follows from the following lemma.

\proclaim{Lemma 3.3}
Let $L_1$ and $L_2$ be two lines in $X=\C^2$ intersecting
transversally at $x=0$. Assume that we have a square
$$\CD
 (\t Y,F) @>{g}>> (Y,y) \\
@VV{\t \pi}V        @VV{\pi }V \\
 (\t X, E) @>{f}>> (X,x)\\ 
\endCD
$$
in which $f$ is the blow up at $x$.
Set $D=L_1+aL_2$ for some $0\le a\le 1$ and assume that
$\pi^* D$ is a Weil divisor.
Then 
$$c_2 (g; {\Olog {\t \pi} {\t X} {f^{-1}D+E}})=0.$$
\endproclaim

\pr
It is sufficient to prove the lemma for one special map
$\pi$.  Let $Y\subset \C ^3$ be a cone over the curve 
$a^n+b^n+c^n=0$ in $\P^2$ and define $\pi \: Y\to X$ by $\pi((a,b,c))=(a^n,b^n)$.
Using Lemma 8.4 one can see that
$${\Olog {\t \pi} {\t X} {f^{-1}D+E}}=\O_{\t Y}\op \O_{\t Y}({\t \pi}^*
(K_{\t X}+f^{-1}D+E)).$$
Hence the second Chern class vanishes (e.g., by Theorem 1.10), Q.E.D.

\medskip
Now we can easily define a global
orbifold Euler number.

\proclaim{Definition 3.4} The \rom{ orbifold Euler number} of the pair 
consisting of a projective surface
$X$ and a $\Q$-divisor $D$ is equal to 
$$\split
e_{\orb}(X,D)=& e_{\top}(X-\Sing X-\Supp D)+ 
\sum_{i} (1-a_i)e_{\top}(D_i-\Sing (X,D))\\
&+\sum_{x\in \Sing (X,D)} e_{\orb}(x;X,D).\\
\endsplit$$
\endproclaim

Note that 
$$e_{\orb}(X,D)=e_{\top}(X)- \sum a_i e_{\top}(D_i-\Sing (X,D))+\sum_{x\in \Sing (X,D)} (e_{\orb}(x;X,D)-1).$$
The formula in Definition 3.4 is more complicated, because it shows how one should think
about the orbifold Euler number: smooth points of $X$ not lying on $D$ should be counted as $1$, smooth points 
of $D_i$ should be counted as $1-a_i$ and the remaining finite set of points should be counted
as appropriate local orbifold Euler numbers. This interpretation 
allows us to define orbifold Euler numbers of $(X,D)$
for any constructible subset of $X$.

\proclaim{Definition 3.5} (see [La3], Definition 2.6)
Let $(X,D)$ be a surface with boundary and $\pi \: Y \to X$ any finite
morphism such that $\pi ^* D$ is a Weil divisor. Then the \rom{ second Chern number}
of $(X,D)$ is defined by 
$$c_2(X,D)= c_2(Y, {\Olog \pi X D})/\deg \pi .$$
\endproclaim
 
One can show that $c_2 (X,D)$ can always be defined and that it 
does not depend on the choice of $\pi$ (see [La3]).

\proclaim {Theorem 3.6}
If $X$ is a projective surface and $D$ a $\Q$-divisor on $X$ then 
$$c_2(X,D)= e_{\orb}(X,D).$$
\endproclaim

\pr
Using Kawamata's covering lemma we can construct a square
$$\CD
 (\t Y,F) @>{g}>> Y \\
@VV{\t \pi}V        @VV{\pi }V \\
 (\t X, E) @>{f}>> X\\ 
\endCD
$$
in which $\pi^*D$ is a Weil divisor and
$f$, $g$ are log resolutions of the pairs $(X, D)$ and
$(Y,\pi^*D)$, respectively.
Then 
$$(g_*({\Olog {\t \pi} {\t X} {(f^{-1}D+E)}}))^{**}=
{\Olog \pi X D}.$$
Hence by Definition-Proposition 1.7
$$c_2(X,D)=c_2({\t X},f^{-1}D+E)-{1\over \deg {\t \pi}}
c_2(g,{\Olog {\t \pi} {\t X} {(f^{-1}D+E)}}).\leqno (3.6.1)$$
Let us compute the first term on the right hand side:
$$\split
&c_2({\t X},f^{-1}D+E)=e_{\orb}({\t X},f^{-1}D+E)=
e_{\top}({\t X}-\Sing ({\t X},f^{-1}D+E))+ \\
&\sum_{i} (1-a_i)e_{\top}(f^{-1}D_i-\Sing ({\t X},f^{-1}D+E))+
\sum_{j} (1-1)e_{\top} (E_j-\Sing ({\t X},f^{-1}D+E))\\
&+\sum_{x\in \Sing ({\t X},f^{-1}D+E)} e_{\orb}(x;{\t X},f^{-1}D+E).
\endsplit
$$
We have
$e_{\top}({\t X}-\Sing ({\t X},f^{-1}D+E))=e_{\top}(X-\Sing (X,D)).$
Since $f^{-1}D_i\cap f^{-1}D_j=\emptyset$ for $i\ne j$, we also have
$e_{\top}(f^{-1}D_i-\Sing ({\t X},f^{-1}D+E))=e_{\top}(D_i-\Sing (X,D)).$
Moreover,
$e_{\orb}(x;{\t X},f^{-1}D+E)=0$
for $x\in \Sing ({\t X}, f^{-1}D+E)$ (by Lemma 3.2).
Therefore 
$$c_2({\t X},f^{-1}D+E)=
e_{\top}(X-\Sing X-\Supp D)+ 
\sum_{i} (1-a_i)e_{\top}(D_i-\Sing (X,D)).
$$
Since the last term  in $(3.6.1)$ is equal to $\sum_{x\in \Sing (X,D)} e_{\orb}(x;X,D)$,
we get the required equality, Q.E.D.

\proclaim{Corollary 3.7}
If $\pi\: Y \to X$ is a  finite proper morphism of normal proper surfaces and
$K_Y+D'=\pi ^* (K_X+D)$ for some boundary $\Q$-divisors then 
$$e_{\orb}(Y,D')=\deg \pi \cdot e_{\orb} (X,D).$$
In particular, if $B$ is the branch locus and $R$ is the ramification locus of $\pi$
then
$$e_{\orb}(Y,R)=\deg \pi \cdot e_{\orb} (X,B).$$
\endproclaim

\pr
The corollary follows immediately from Theorems 1.6, 3.6 and the logarithmic ramification formula:
$${\Olog \pi X D}= {\olog Y {D'}}.$$
The last equality holds on each finite covering of $Y$ on which both sheaves are defined, Q.E.D.

\medskip

\noindent
{\sl 3.8. Remarks.}

\item{1.} Corollary 3.7, together with explicit formulas for local orbifold 
Euler numbers (see Theorems 6.3, 8.3, 8.7, 9.4.2), allows us to  compute
changes of orbifold Euler numbers for large classes of finite morphisms. 
In special cases we get results of Holzapfel [Hz], who considered only locally simple cases,
e.g., when the branch locus is a normal crossing divisor (Proposition 9.3.1 is sufficient
to deal with this case).

\item{2.} If $\pi \: Y\to X$ is a Galois covering then all the ramification indices over a fixed
component of the branch locus, say $B_i$, are equal (to $b_i$). Then
$$e_{\orb}(Y)=\deg \pi \cdot e_{\orb} (X,\sum \left(1-{1\over b_i}\right)B _i).$$ 

\medskip

\heading 4. Logarithmic Ramification Formula and a vanishing theorem
for log canonical pairs \endheading 

In 4.1--4.7 we study a conjectural characterization of log canonical pairs
in terms of symmetric powers of logarithmic $1$-forms (it is also possible
to extend it to $q$-forms). However, the main aim of this section
is a proof of the Logarithmic Ramification Formula (Theorem 4.9)
and a Bogomolov--Sommese type vanishing theorem (Theorem 4.11)
for log canonical pairs. These results are used in the next section
in the proof of inequalities between logarithmic orbifold Euler numbers.

\proclaim{Definition 4.1} 
Let $(X,D)$ be a pair where $X$ is a normal variety and $D$ is a boundary $\Q$-divisor.
Consider a square
$$\CD
 {\t Y} @>{g}>> Y \\
@VV{\t \pi}V        @VV{\pi }V \\
 (\t X ,E) @>{f}>> X\\ 
\endCD
$$
such that ${\t \pi}^*(f^{-1}D+E)$ is a Weil divisor and $\t X$ is a log resolution of
$(X,D)$.
The pair $(X,D)$ is called \rom{1-log canonical} (or \rom{1-lc}) if the 
sheaves 
$g_*\S ^n {\Olog {\t \pi} {\t X} {f^{-1}D+E}}$ are reflexive for all $n$.
If $Z\subset X$ is a closed subvariety then $(X,D)$ is called \rom{1-log canonical 
along} $Z$ if the pair $(X,D)$ is 1-lc in an open neighbourhood of $Z$.
\endproclaim

A smooth log pair is 1-lc by the logarithmic ramification formula.
The condition $(X,D)$ log canonical is equivalent to the similar condition with $1$-forms replaced by $m$-forms, 
where $m=\dim X$. It is not clear from the definition if it depends on the choice of 
$\pi \: Y\to X$, so we fix $\pi$  in  further considerations. We will show that
the definition does not depend on the choice of resolution $f$.

\proclaim{Lemma 4.2}
Set $\E ={\Olog \pi X D}$ and $\F={\Olog {\t \pi} {\t X} {f^{-1}D+E}}$.
The following conditions are equivalent:
\item{(1)} $g_*\S^n\F$ is reflexive,
\item{(2)} $g_* \S^n\F=\S ^n\E$,
\item{(3)} $(g^*\S^n\E )^{**}\hookrightarrow \S^n \F$,
\item{(4)} $(g^*(\S^n \E ))^{**}\hookrightarrow  {\t \pi }^*(S^n{\olog {\t X} {\lceil f^{-1}D\rceil+E}})$,
\item{(5)} There exists an analytic subset $S$ of points on the $g$-exceptional divisor $F$
which has codimension at least $2$ in $\t Y$ and such that for every $\omega \in \S ^n\E$  
the germ of $g^* \omega$ at any point $P\in F-S$ belongs to $({\t \pi }^*S^n{\olog {\t X} {E}})_P$.
\endproclaim

\pr 
The equivalence of (1) and (2) is clear. Since $g_*((g^*\S^n\E )^{**})$ is reflexive we get the equivalence
of (2) and (3). Obviously, (3) implies (4) and (4) implies (5). To prove that (3) implies (5)
note that $(g^*\S^n\E )^{**}$ is an $\O _{\t Y}$-submodule of $n$-th symmetric power of meromorphic
$1$-forms on $\t Y$ generated by $g^* \omega$ for $\omega \in \S ^n \E$.
Hence it is sufficient to prove that the germ of $g^* \omega $ 
belongs to $(\S^n\F)_{P, {\t Y}}$  for every point $P\in {\t Y}$. 
But this is true for $P\in {\t Y}-S$, so also everywhere. Q.E.D.

\medskip

\re 

(1) For $n=1$ conditions (3) and (4) can be thought of as a generalization of the logarithmic ramification formula
to the singular case.
 
(2) Condition (4) in the above lemma  and the logarithmic ramification formula say
that we can check if $(X,D)$ is 1-lc using only one fixed log resolution.

\medskip

\noindent
{\sl 4.3.} If we have a diagram
$$\CD
{\t Z} @>{h}>> Z \\
@VV{\t \varphi}V        @VV{\varphi }V \\
 {\t Y} @>{g}>> Y \\
@VV{\t \pi}V        @VV{\pi }V \\
 (\t X ,E) @>{f}>> X\\ 
\endCD
$$
and $h_*S ^n {\Olog {{(\t \pi}{\t \varphi})} {\t X} {f^{-1}D+E}}$ is reflexive then
$g_*S ^n {\Olog {\t \pi} {\t X} {f^{-1}D+E}}$ is also reflexive.
This follows from the fact that $\O _{\t X}$ is a direct summand of ${\t \pi}_*\O _{\t Y}$
(see, e.g, [KM], Proposition 5.7).

Therefore to check if  $(X,D)$ satisfies (1) of Lemma 4.2 for all squares it is sufficient 
to check the equivalent conditions of Lemma 4.2 for all maps $Y\to X$ which  factor through a given map. 
In particular,  in the algebraic case we can use Kawamata's covering trick to assume
that ${\Olog {{\t \pi}} {\t X} {f^{-1}D+E}}$ is a locally free sheaf.
Moreover, we have the following lemma:

\proclaim{Lemma 4.4}
Let  $\pi \: Y\to X$ be a finite covering and $K_Y+{D'}=\pi ^* (K_X+D)$
for some boundary divisors $D$ and $D'$. Then $(X,D)$ is 1-lc for all squares
if and only if $(Y,D')$ is 1-lc for all squares.
\endproclaim

\proclaim{Proposition 4.5}
Let $(X,D)$ be a 1-lc pair and let $\pi \:Y \to X$ be a finite map such that $\pi ^*D$ is a Weil divisor.
If a rank $1$ reflexive sheaf $\L$ is contained in $\OLog \pi q X D$ then 
$\kappa (Y, \L)\le q$.
\endproclaim

\pr
Note $\F={\t \pi}^* {\oLog q {\t X} {\lceil f^{-1}D \rceil +E} }$ is locally free
and by Lemma 4.2  it contains $(g^*\L )^{**}$. Take a rank 1 reflexive subsheaf $\M$ of $\F$
containing  $(g^*\L )^{**}$ and such that $\q= \F/\M$ is torsion free. Passing to the 
resolution of singularities of $\t Y$ dominating the Nash blow up of $\q$ we can assume
that $\M$ and $\q$ are locally free. It follows that $\M^{\ot n}$ is a subline bundle
of $S^n\F$ with a locally free quotient $\q _n$. 
Now we can use the diagram  
$$\CD
 0 @>>> H^0({\t Y},\M^{\ot n}) @>>> H^0({\t Y}, S^n \F) @>>> H^0({\t Y}, \q _n)\\
@. @VVV      @VVV  @VVV \\
0 @>>> H^0({\t Y}-F,\M^{\ot n}) @>>> H^0({\t Y}-F, S^n \F) @>>> H^0({\t Y}-F, \q _n)\\ 
\endCD
$$
to deduce that $H^0({\t Y},\M^{\ot n})=H^0({\t Y}-F,\M^{\ot n})=H^0({Y},\L ^{\ot n})$.
Since $\kappa ({\t Y}, \M )\le q$ by Theorem 4.10, we also get $\kappa ({Y}, \L  )\le q$, Q.E.D.

\proclaim{Proposition 4.6}
Let $(X,D)$ be a normal surface with boundary.
If $(X,D)$ is 1-log canonical then it is log canonical.
\endproclaim

\pr
The problem is local so we can work on a normal surface germ $(X,x)$.

Set $\F={\Olog \pi X D}$ and $\G={\Olog {\t \pi} {\t X} {f^{-1}D+E}}$. By assumption
$g_*(S^n\G)=\S ^n \F$. Therefore we have an inclusion
$$(\S ^n \F )^p=(g^* (\S ^n \F ))^{**}\hookrightarrow S^n\G .\leqno (4.6.1)$$
Since $g_*c_1(S^n \G )=c_1(\S ^n \F )$ the divisor 
$D_n=c_1(S^n \G )-g^*c_1(\S ^n \F )={n(n+1) \over 2} (c_1 \G -g^* c_1\F )$ is 
supported on the exceptional set of $g$.
Now $(4.6.1)$ implies that $D_n\ge  
c_1((\S^n \F)^p)-g^*c_1(\S ^n \F )= c_1(g, (\S ^n \F )^p)$. But by Proposition 4.6, [La1],
divisors
$c_1(g, (\S ^n \F )^p)$ have uniformly bounded coefficients. Hence $c_1\G \ge g^*c_1\F$ and
$K_{\t X}+f^{-1}D+E \ge f^* (K_X+D)$, so $(X,D)$ is lc, Q.E.D.

\medskip
By Corollary 7.5 the converse of Proposition 4.6 holds  if the pair $(X,D)$
has at most fractional boundary.

\proclaim{Problem 4.7}
Is it true that $(X,D)$ is log canonical if and only if it is 1-log canonical for all squares? 
\endproclaim

By Theorem 4.9 if there exists $\pi \: Y\to X$ such that $\Olog \pi X D$ is locally free
and  $(X,D)$ is log canonical then it is 1-log canonical. In  particular,
quotients of smooth pairs are 1-lc.

\proclaim{Lemma 4.8}
Let $\pi \:Y\to X$ be a finite map of smooth varieties branched only along a smooth divisor $E$.
Let $i\: Y-\pi ^{-1}(E)\hookrightarrow Y$ denote an embedding. 
If $\L$ is an $\O _Y$-submodule of $i_*\pi^* \Omega ^q_X$ and $d(\L)\subset \pi^*{\oLog {{q+1}} X E}$
then $\L\subset \pi^*{\oLog q X E}$.
\endproclaim

\pr
The statement is local so we can work on germs only. For simplicity let us assume that $q=1$. In general
the proof goes along the same lines.

Assume that there exists
a meromorphic $1$-form $\omega \in \L$ which does not belong to $\pi^*{\oLog q X E}$.
Let $z=0$ be an equation of $E$ in $X$. By assumption there exists a positive integer $k$ such that 
$\pi^*(z^k)\omega \in \pi^* \Omega _X$. Let us take minimal such $k$.
If $k=1$ then $\omega \in \pi^*{\oLog q X E}$, a contradiction. Therefore $k\ge 2$.  

By assumption $d(\pi ^* (f)\omega )\in \pi^*{\oLog {{q+1}} X E}$ for any $f\in \O _X$.
Therefore 
$$d(\pi^*(f))\wedge (\pi^*(z)\omega )= \pi^*(z)(d(\pi^*(f)\omega)-\pi^*(f)d\omega ) \in 
{\pi ^*\Omega^{q+1}_X }\leqno (4.7.1)$$
for any $f\in \O_X$. Let us write $\pi^*(z)\omega$ using a local coordinate system $x_1,\dots, x_n$
on $X$:
$$\pi^*(z)\omega=f_1\pi ^*(dx_1)+\dots +f_n\pi^* (dx_n).$$
By (4.7.1) the form $d(\pi^*(x_i))\wedge (\pi^*(z)\omega)= 
\sum _{j\ne i}f_j \pi^*(dx_i)\wedge \pi^*(dx_j)$ belongs to ${\pi ^*\Omega^{{q+1}}_X }$.
Therefore $f_i$ belongs to $\O _Y$ for $i=1,\dots , n$ and $\pi^*(z)\omega \in \pi^* \Omega _X$, 
a contradiction, Q.E.D. 

\medskip
Clearly, the lemma is false on the level of elements, e.g., if $\pi$ is the identity then
$dz\over z^n$ is a closed $1$-form which does not belong to ${\oLog 1 X {(z=0)}}$.

\proclaim{Theorem 4.9} (Logarithmic Ramification Formula for log canonical pairs)
Let us assume that $(X,D)$ is log canonical and consider the square
$$\CD
 {\t Y} @>{g}>> Y \\
@VV{\t \pi}V        @VV{\pi }V \\
 (\t X ,E) @>{f}>> X\\ 
\endCD
$$
such that ${\t \pi}^*(f^{-1}D+E)$ is a Weil divisor and $f$ is a log resolution of
$(X,D)$. Then 
$$(g^*{\OLog \pi q X D})^{**}\hookrightarrow {\OLog {\t \pi} q {\t X} {f^{-1}D+E}}$$
and $g_*{\OLog {\t \pi} q {\t X} {f^{-1}D+E}}$ is reflexive for any $1\le q\le n=\dim X$.
\endproclaim

\pr
Set $\E ^q={\OLog {\pi} q X D}$ and $\F ^q={\OLog {\t \pi} q {\t X} {f^{-1}D+E}}$.
The sheaves $(g^*\E ^q)^{**}$ and $\F ^q$ are subsheaves of the sheaf $\M _{\t Y}^q$ 
of meromorphic $q$-forms on $\t Y$.
Therefore to show that $g_*\F ^q$ is reflexive it is sufficient to show that there exists an
analytic subset $S$ of $\t Y$ of codimension $\ge 2$ such that for any $y\in {\t Y}-S$
the 
inclusions $(g^*\E ^q)^{**}\subset \M ^q_{\t Y}\ot \Omega _{\t Y}$ and $\F ^q\subset \M ^q_{\t Y}\ot \Omega _{\t Y}$ 
induce an inclusion of stalks
$$(g^*\E ^q)^{**}_{y, {\t Y}}\hookrightarrow (\F ^q)_{y, {\t Y}}.\leqno (*)_q$$
If $y\in {\t Y}-F$ then we have equality in $(*)_q$. There exists a codimension $\ge 2$ subset $S$ of $\t Y$
which contains $F\cap {\t \pi}^{-1}D$ and such that in a neighbourhood of any point of $E$
the map ${\t \pi} \: {\t Y -S}\to {\t \pi}({\t Y}-S)$  is  branched only along smooth points of $E$.

Note that $(*)_n$ follows from the assumption that $(X,D)$ is log canonical.
Therefore for any point $y\in F-S$ Lemma 4.8 implies $(*)_q$ by induction on $k=n-q$, Q.E.D.

\medskip
A special case of Theorem 4.9, when $D$ is a reduced divisor
on a normal surface, was proved in [La3], Theorem 4.2.

\proclaim{Theorem 4.10} (Bogomolov--Sommese vanishing theorem; [EV], Corollary 6.9)
Let $X$ be a smooth projective manifold and $D$ a normal crossing divisor on $X$. Let $\pi \: Y \to X$
be a generically finite map from a normal variety $Y$. If a line bundle $\L$ is contained
in $\pi ^* {\oLog q X D}$ then $\kappa (Y, \L) \le q.$
\endproclaim

\pr
Let $g\: (Z, B)\to (Y,(\pi ^* D)_{red})$ be a log resolution. By the logarithmic ramification
formula $(\pi g)^*{\oLog q X D}\hookrightarrow {\oLog q Z B}$. Hence $g^* \L \subset {\oLog q Z B}$
and the assertion follows from [EV], Corollary 6.9, Q.E.D.

\medskip
The following theorem is a Bogomolov--Sommese type vanishing theorem for log canonical $\Q$-divisors.

\proclaim{Theorem 4.11} 
Let $X$ be a normal projective variety and assume that $(X,D)$ is log canonical.
Let $\pi \:Y \to X$ be a finite map such that $\pi ^*D$ is a Weil divisor.
If a rank $1$ reflexive sheaf $\L$ is contained in $\OLog \pi q X D$ then 
$\kappa (Y, \L)\le q$.
\endproclaim

\pr
Assume that $\kappa (Y, \L )>q$. Write $\L= \O _Y(L)$ for some Weil divisor $L$.  
Let $H$ be any very ample Cartier divisor on $Y$.
By Kodaira's lemma ([KMM, Lemma 0-3-3]) there exists an integer $n$ such that
the linear system $|nL -H|$ contains a divisor $D$. By the covering lemma there exists
a finite map $p\:Z\to Y$ such that $p^*({1\over n }D)$ is an effective Weil divisor.
Now we can find another covering $q\: T\to Z$ and a Cartier divisor $H'$ on $T$ 
such that $q^*(p^*H)$ is linearly equivalent to $nH'$
(see [KM], Proposition 2.67). Hence $\O _T(H')\subset 
((pq)^*\L)^{**}\subset {\OLog {(\pi pq)} q X D}$, where $ H'$ is an ample Cartier divisor.
Now passing to a log resolution of $(X,D)$, forming a square and using Theorem 4.9
we get a contradiction with Theorem 4.10, Q.E.D.

\medskip

\re 
Theorem 4.11 is similar to Proposition 4.5 but it is not equivalent to it. Proposition 4.5
uses the  conditions of Lemma 4.2 for one fixed square whereas Theorem 4.11 uses these conditions for
$n=1$ and for all squares.

\medskip

\heading 5. Inequalities between logarithmic orbifold Euler numbers\endheading

In this section we prove Theorems 0.1 and 0.3 from the Introduction.

\proclaim{Theorem 5.1}
Let $(X,D)$ be a log canonical surface pair  and let $\pi \: Y\to X$ be a finite covering such that
$\pi^*D$ is a Weil divisor. Let $\F \subset {\Olog \pi X D}$ be a rank 
2 reflexive subsheaf such that $c_1\F$ is pseudoeffective. Let $c_1\F=P+N$
be the Zariski decomposition. Then
$$3c_2\F+{1\over 4} N^2\ge (c_1\F)^2.$$
\endproclaim

The proof of this theorem is the same as the proof of Theorem 0.1 in [La3].
The only difference is that we need to replace Corollary 4.4, [La3] with the more general
Theorem 4.11.  As an immediate corollary of Theorems 3.6 and 5.1 we get the following generalisation
of Theorem 0.1.

\proclaim{Corollary 5.2}
If $(X,D)$ is lc and $K_X+D=P+N$ is pseudoeffective then
$$3e_{\orb}(X,D)\ge 3e_{\orb}(X,D)+{1\over 4}N^2\ge (K_X+D)^2.$$
\endproclaim

Note that if $3e_{\orb}(X,D)= (K_X+D)^2$ then $N^2=0$. Hence $N=0$ and $K_X+D$ is nef.

\medskip
Theorem 0.3 is obtained from the above result by using the formula 
$e_{\top} (D_i)=2-2g(D_i)-\sum _{P\in \Sing D_i}(r_P-1)$ and Corollary 7.8.

\medskip

\heading 6. Orbifold Euler numbers for non lc pairs and logarithmic 
cohomology of open surfaces \endheading

In this section we compute Chern classes of sheaves of logarithmic 1-forms with poles along reduced divisors
on smooth varieties and apply them to study the de Rham morphism for open surfaces.

\medskip

Let $c(\E)$ denote the total Chern class of the sheaf $\E$. This can be defined for any sheaf on a smooth 
projective variety.

\proclaim{Proposition 6.1}
Let $D$ be a hypersurface in a smooth projective variety $X$. Then
$$c(({\hlog X D})^*)={c(T_X)\cdot c(\O_Y(D))\over c(\O_X(D))},$$
where $Y$ denotes the singular subscheme of $D$ (i.e., the subscheme of $D$ defined by 
the partial derivatives of a section defining $D$).
\endproclaim

\pr
Let $\J ^1_XD$ be the bundle of 1-jets of $\O _X(D)$. It  is defined by the
divisor class $[D]\in H^1(\Omega _X)\simeq \Ext ^1(\O _X(D), \Omega _X(D))$.
The divisor $D$ determines a section $\O _X \to \J ^1_X D$. Let $\E$ be the cokernel
of this section after twisting by $\O _X(D)$. Then $({\hlog X D})^*=\E ^*$ fits into
the following exact sequence
$$0\raa \E^* \ra (\J ^1_X D)^*\ot \O _X(D)\ra \O _X(D)\ra \O _Y (D)\raa 0 . $$
Now the result follows from the sequence
$$0\raa \O_X \ra (\J ^1_X D)^* \ot \O _X(D) \ra T _X\raa 0\hbox{, Q.E.D.}$$

\medskip
Let $C$ be a curve on a smooth surface germ $(X,x)$. Let $z_1, z_2$ be local
coordinates at $x$. If $C$ is given by $f(z_1,z_2)=0$ then we define the {\sl Tjurina number} of $C$ at $x$
by $\tau _x(C)=\dim _{\C} \O _{X , x}/(f, {{\partial f}\over \partial z_1},
{{\partial f}\over \partial z_2})$ and the {\sl Milnor number} of $C$ at $x$ by
$\mu _x(C)=\dim _{\C} \O _{X,x}/( {{\partial f}\over \partial z_1},
{{\partial f}\over \partial z_2})$.

\proclaim{Corollary 6.2}
Let $D$ be a reduced divisor on a smooth surface $X$. Then 
$$c_2({\hlog X D})=e_{\top }(X-D)+\sum _{x\in \Sing D}(\mu _x(D)-\tau _x(D)).$$ 
In particular, $c_2({\hlog X D})=e_{\top }(X-D)$  if and only if $D$ has locally
weighted homogeneous singularities.
\endproclaim

\pr
By Proposition 6.1, $c_2({\hlog X D})=c_2(X)+(K_X+D)D-\sum _{x\in \Sing D}\tau _x(D)$.
Since $e_{top} (X-D)=e_{\top}(X)-e_{\top}(D)=c_2(X)+(K_X+D)D-\sum _{x\in \Sing D}\mu _x(D)$, 
we get the first part of the corollary.
The second part follows from the first  and from
Saito's theorem (see [Sai1]) saying that $\tau _x=\mu _x$ if and only if $D$ has a
weighted homogeneous singularity at $x$, Q.E.D.

\proclaim{Corollary 6.3}
Let $C:=(f=0)$ be a curve in $(\C ^2, 0)$. Then
$$e_{\orb}(0; \C^2, C) =\mu _0 (C)-\tau _0(C)\ge 0.$$
\endproclaim

\pr
We give a proof by globalisation, since it is much easier than a
local proof. By the finite determinacy property we can find 
a smooth projective surface $X$ and a curve $D$
on $X$ singular at only one point $x$ such that the germs $(X,D,x)$ and
$(\C ^2, C, 0)$ are of the same analytic type. Let $f\: \t X \to X$
be a log resolution of $(X, D)$ at $x$. Then 
$$\split
e_{\orb}(x; X, D)&=-c_2(x, {\olog {\t X} {f^{-1}D+E}})=
c_2({\hlog X D})-c_2({\olog {\t X} {f^{-1}D+E}})\\
&=c_2({\hlog X D})-e_{top} (X-D)\\
\endsplit$$
and the assertion follows from Corollary 6.2, Q.E.D.

\medskip

Let $D$ be a divisor on a smooth complex manifold $X$. Let $U$ be the complement of $D$ in $X$
and $j\:U\hookrightarrow X$ the inclusion. We say that {\sl the logarithmic comparison theorem holds for} 
$D$ if the de Rham morphism gives rise to an isomorphism (in the derived category)
$${\hLog \bullet X D} {\ra^{\simeq}}{\Bbb R}j_* (\C _U)$$
(see [CCMN]).

\proclaim{Corollary 6.4} (Theorem 1.3, [CCMN])
Let $C$ be a plane curve. If the logarithmic comparison theorem holds for $C$
then all the singularities of $C$ are locally weighted homogeneous.  
\endproclaim

\pr
The problem is local so we can fix a point $P$ on $C$.
By the finite determinacy property  we can find a smooth projective surface $X$ and a curve $D$
on $X$ which is singular at only one point $x$ such that the germs $(X,D,x)$ and
$(\C ^2, C, P)$ have the same analytic type.

By the assumption, if $U=X-D$ and $j\: U \hookrightarrow X$ is the inclusion,
then the de Rham morphism gives rise to an isomorphism 
$${\hLog \bullet X D} {\ra^{\simeq}}{\Bbb R}j_* (\C _U).$$
Therefore the induced map on the hypercohomology is an isomorphism
$${\Bbb H} ^*(X,{\hLog \bullet X D}) \ra^\simeq {\Bbb H}^* (X,{\Bbb R}j_* (\C _U)).$$
By standard arguments (the Poincar\'e lemma and spectral sequences) 
${\Bbb H} ^* (X,{\Bbb R}j_* (\C _U))\simeq H^*(U,\C ).$
On the other hand we have the spectral sequence 
$$E^{pq}_2=H^q (H^p (X,{\hLog \bullet X D}))$$
abutting to ${\Bbb H} ^*(X,{\hLog \bullet X D})\simeq H^*(U,\C )$. 
Comparing Euler characteristics for complexes we get
$$\split
e_{\top}(U)&=\sum _p (-1)^ph^p (U, \C)=
\sum _{p,q}(-1)^{p+q}\dim E^{pq}_2=\sum _{p,q}(-1)^{p+q}h^{p}(X,{\hLog q X D})\\
&=\sum _q(-1)^q \chi (X,{\hLog q X D}).\\
\endsplit$$
By the Riemann--Roch theorem the last expression is equal to 
$c_2({\hlog X D})$. Hence the corollary follows from Corollary 6.2, Q.E.D.

\medskip

\heading 7. Properties of orbifold Euler numbers\endheading

The main aim of this section is to prove properties (0.2.1)--(0.2.3) of local orbifold Euler numbers
stated in the Introduction.

\proclaim{Lemma 7.1}
Let  $\pi \: (Y,y)\to (X,x)$ be a finite proper map of normal surface germs. 
If $K_Y+{D'}=\pi ^* (K_X+D)$
for some boundary divisors $D$ and $D'$ then
$$e_{\orb}(y;Y,D')=\deg \pi \cdot e_{\orb}(x;X,D) .$$
\endproclaim

\re
In particular, if $\pi$ is \'etale then we can reduce the computation of
$e_{\orb}(x;X,D)$ to that of $e_{\orb}(y;Y,\pi ^*D)$.

One can also reduce the study of log canonical pairs to log terminal pairs using
log crepant morphisms:

\proclaim{Lemma 7.2} (cf.~Lemma 4.6, [Me2])
Let $(X,D)$ be a normal projective surface with log canonical singularities.
Then there exists a birational morphism $\psi : ({\hat X}, {\hat D})\to (X, D)$
such that ${\hat D}$ is the sum of the strict transform of $D$ and the exceptional
curves with coefficient $1$, $K_{\hat X}+{\hat D}= \psi ^* (K_X+D)$,
$c_2(X,D)\le c_2({\hat X}, {\hat D})$ and  $({\hat X}, {\hat D})$ has log terminal
singularities.
\endproclaim

The lemma follows from the following local statement:

\proclaim{Lemma 7.3}
Let $(X,D,x)$ be a germ of a log canonical surface pair.
Assume that $(X,D)$ is not log terminal at $x$. Then there exists 
a square
$$\CD
 {\hat Y} @>{\hat \psi}>> Y \\
@VV{\hat \pi}V        @VV{\pi }V \\
 (\hat X ,F) @>{\psi}>> (X,x)\\ 
\endCD
$$
such that  $c_1({\hat \psi}, {\Olog {\hat \pi} {\hat X} {\hat D}})=0$,
$c_2({\hat \psi}, {\Olog {\hat \pi} {\hat X} {\hat D}})\ge 0$
and $({\hat X}, {\hat D})$ is log terminal, where ${\hat D}=\psi^{-1} D+F$.
\endproclaim

\pr
Let us take the minimal log resolution $f\: ({\t X}, E)\to (X,x)$ of $(X,D)$
and write 
$$K_{\t X}+f^{-1}D+E- f^*(K_X+D)=\sum a_i E_i$$
for some $a_i\ge 0$. Let $\varphi \: {\t X}\to {\hat X}$
be the contraction of all the exceptional curves with $a_i\ne 0$ and let $\psi \: {\hat X} \to X$
be the induced morphism. In some cases we do not contract any curve, e.g., if $X$
is a cone over a plane conic and $D$ is a sum of two generators with coefficient $1$.
Note that $F=\varphi _* E$ and
$$K_{\hat X}+\hat D=\varphi _*(K_{\t X}+f^{-1}D+E)=\varphi _*(f^*(K_X+D))=\psi ^*(K_X+D).$$
Therefore, since $(X,D)$ is lc  $({\hat X}, {\hat D})$ is log terminal.  
Finally,
$$c_2({\hat \psi}, {\Olog {\hat \pi} {\hat X} {\hat D}})\ge {1\over 4}c_1^2
({\hat \psi}, {\Olog {\hat \pi} {\hat X} {\hat D}})=0$$
by Definition 1.3, Q.E.D.
\medskip
Note that $({\hat X}, {\hat D})$ is not klt at the points of $F$ and 
${\hat X}$ has at most cyclic quotient singularities. In some cases we can prove the equality 
$c_2({\hat \psi}, {\Olog {\hat \pi} {\hat X} {\hat D}})=0$ in Lemma 7.3.

\proclaim{Proposition 7.4}
Under the assumptions of Lemma 7.3 assume that the singularities of $({\hat X}, {\hat D})$ are locally 
quotients of smooth log pairs. Then $e_{\orb} (x;X,D)=0$ and $(X,D)$ is 1-lc (for all squares).
\endproclaim

\pr
We can use coverings and simple base change to construct a square
$$\CD
{\hat Y} @>{\hat \psi}>> Y\\
 @VV{\hat \pi}V  @VV{\pi }V \\
 {\hat X} @>{\psi}>> X\\ 
\endCD
$$
such that ${\hat \pi}^* D$ is a Weil divisor and $\F= {\Olog {\hat \pi} {\hat X} {\hat D}}$ 
is a locally free sheaf. Set $\G= {\Olog {\pi} {X} {D}}$. 
Let $\nu \: Z\to {\hat Y}$ be a resolution of singularities of $\hat Y$ and set ${\t \psi}={\hat \psi}\nu$.
We can use the residue map to show that if $F_i$ is a smooth ${\t \psi}$-exceptional divisor
then there exists a surjection 
$$\nu ^*\F |_{F_i}\to \O _{F_i}.$$
The kernel of this map is a degree $0$ line bundle on $F_i$. By [Wa2], Proposition 3.16
$$\chi ({\t \psi},S ^n \nu ^* \F)=\chi ( {\hat \psi},S ^n \F)=\bo (n^2).$$
Together with $c_1({\t \psi},\nu ^* \F )=0$ this gives  
$c_2({\hat \psi}, \F)=c_2({\t \psi},\nu ^* \F )=0$. Now the first statement follows from
Lemma 3.2 and Lemma 7.1.
To prove the second part note that $S^n \nu ^*\F |_{F_i}$ has a filtration by line bundles of degree $0$ and
for any negative degree line bundle
$\L$ on $F_i$ we have $H^0(S^n (\nu^*\F )|_{F_i}\ot \L)=0$. 
Hence we can use the formal function theorem to show that $H^1_F(S^n(\nu ^*\F))=0$, where
$F$ is  the whole ${\t \psi}$-exceptional set with a reduced structure (see [Wa1], Proposition 2.4). 
Therefore ${\hat \psi}_*(S^n\F)={\t \psi}_*(S^n\nu^*\F)=\S ^n \G$ and $(X,D)$ is 1-lc
since $({\hat X}, {\hat D})$ is 1-lc (see Lemma 4.4), Q.E.D.

\proclaim{Corollary 7.5} (cf. Theorem 4.2, [La3])
Let $(X,D)$ be a log canonical surface with at most fractional boundary. Then $(X,D)$ is 1-lc for all squares. 
Moreover, if $(X,D)$ is not lt at $x$ then $e_{\orb}(x;X,D)=0$. 
\endproclaim

\proclaim{Theorem 7.6}
Let $(X,D,x)$ be an affine germ of a log canonical surface with boundary. If $D$ contains a prime divisor $L$ 
with coefficient $1$ then $e_{\orb} (x;X,D)=0.$
\endproclaim

\pr
The proof is by induction on the number $k(X,D)$ of curves in the minimal log resolution of $(X,D)$.

If $k(X,D)=0$ then the theorem holds. If $(X,D)$ is log canonical at $x$ then there exists an effective  
$\Q$-divisor $D'$ such that $(X,D+D')$ is log canonical
but not log  terminal at $x$ and $k(X,D+D')=k(X,D)$ (one can construct $D'$, e.g., by increasing $\alpha$ 
in $(X, L+\alpha (D-L))$).

Let us apply Lemma 7.3 to the pair $(X, D+D')$. Note that 
$c_2({\hat \psi}, {\Olog {\hat \pi} {\hat X} {\hat D}})\ge 0$,
where ${\hat D}= \psi ^{-1}(D+D')+F$.
Hence by the induction assumption and Definition--Proposition 1.5 
$$e_{\orb}(x;X,D+D')=\sum {} e_{\orb}(z; {\hat X}, {\hat D})-
c_2({\hat \psi}, {\Olog {\hat \pi} {\hat X} {\hat D}})/\deg {\hat \pi}\le 0.\leqno (7.6.1)$$
Taking a finite covering of $X$ we can assume that $X=\C^2$ and 
$L$ is a line passing through $x=0$.

\medskip
{\sl Claim 7.6.2.} 
Let $f\: {\t X}\to X$ be any resolution of $(X,x)=(\C^2, 0)$. Then we have a short exact sequence
$$0\raa \O_{\t X}\ra^{\omega} {\olog {\t X} {f^{-1}L+E}} \ra \O (K_{\t X}+{f^{-1}L})\raa 0.$$
Moreover, in a neighbourhood of every point $P$ of $E+{f^{-1}L}$ there exists a function $\alpha$
such that  $E+{f^{-1}L}$ is given by $\alpha =0$ and $\omega={d\alpha \over \alpha}$. 

\medskip
\pr The proof is by induction on the number of irreducible components of the exceptional set of $f$.
If $f$ is the blow up  of $X$ at $x$ then we can use Lemma 8.2. 

Let $g\: ({\t Y}, F)\to ({\t X}, E)$ be the blow up of a point $P\in E$. Set $h=gf$. 
If $P$ is a double point of $E$ or the point of intersection of $E$ and $f^{-1}L$ then
$${\olog {\t Y} {h^{-1}L+F}}=g^*{\olog {\t X} {f^{-1}L+E}}$$
and the claim is clear.
Otherwise, let $z_1, z_2$ be local coordinates around $P$, where $E=\{z_1=0\}$ and $\omega ={d z_1\over z_1}$.
On the part of $\t Y$ with coordinates $z_2$ and $t={z_1\over z_2}$ the sheaf
${\olog {\t Y} {h^{-1}L+F}}$ is generated by ${dz_2\over z_2}$ and $g^*\omega ={dt\over t}+{dz_2\over z_2}$.
On the other part of $\t Y$ the claim is clear, Q.E.D.

\medskip
Let $f\: {\t X}\to X$ be any log  resolution of $(X,D)$.
We can find a square
$$\CD
{\t Y} @>{g}>> Y\\
 @VV{\t \pi}V  @VV{\pi }V \\
 {\t X} @>{f}>> X\\ 
\endCD
$$
such that $\E = {\Olog {\t \pi} {\t X} {f^{-1} (D+D')+E}}$ and ${\Olog {\t \pi} {\t X} {f^{-1} D+E}}$
are locally free. Then by the claim $\omega$ gives rise to a short exact sequence
$$0\raa \O_{\t Y}\ra^{{\t \pi}^*\omega} \E \ra \O ({\t\pi}^*(K_{\t X}+{f^{-1}(D+D')}+E))\raa 0.$$
Since $(X, D+D')$ is lc $g_*\E={\Olog {\pi} {X} {D+D'}}$ by Theorem 4.9.
Since $g_*{\O _{\t Y}}=\O _Y$ we have the induced map 
$\O _Y \to {\Olog {\pi} {X} {D+D'}}$. Then in the notation of [La4] we have
$$e_{\orb}(x;X,D+D')={1\over \deg \pi}\delta _y (\O _Y \to {\Olog {\pi} {X} {D+D'}}).$$
Similarly
$$e_{\orb}(x;X,D)={1\over \deg \pi}\delta _y (\O _Y \to {\Olog {\pi} {X} {D}})$$
and by the construction ${\Olog {\pi} {X} {D}}={\Olog {\pi} {X} {D+D'}}[-\pi ^*D']$
is obtained as the elementary transformation of ${\Olog {\pi} {X} {D+D'}}$ with respect
to $\pi ^*D'$. Moreover, ${\Olog {\pi} {X} {L}}$ is obtained as the elementary transformation of 
${\Olog {\pi} {X} {D}}$ with respect to $\pi ^*(D-L)$. Hence by Theorem 2.6, [La4]
$$\delta _y (\O _Y \to {\Olog {\pi} {X} {D+D'}})\ge \delta _y (\O _Y \to {\Olog {\pi} {X} {D}})\ge
\delta _y (\O _Y \to {\pi}^*{\olog  {X} {L}}),$$
that is,
$$e_{\orb}(x;X,D+D')\ge e_{\orb}(x;X,D)\ge e_{\orb}(x;X,L)=0.$$
But $e_{\orb}(x;X,D+D')\le 0$ by (7.6.1) and hence $e_{\orb}(x;X,D)=0$, Q.E.D.

\proclaim{Corollary 7.7}
Let $D$ and $B$ be effective $\Q$-divisors on a normal surface germ $(X,x)$. Assume
that $(X, D+B)$ is log canonical but not log terminal. Let $f\: \t X \to X$ be a log
resolution of $(X, D+B)$. Then
$$e_{\orb}(x; X,D)\le -{1\over 4}(c_1 (f, f^{-1}B))^2.$$ 
\endproclaim

\pr
Let us apply Lemma 7.3 to the pair $(X,D+B)$ and let $\varphi \: {\t X}\to {\hat X}$ 
be such that $c_1({\psi}, K_{\hat X}+{\hat D}+{\psi ^{-1}B})=0$, where ${\hat D}={\psi ^{-1}}D+F$
and $F$ is the exceptional divisor.
Let $\psi \: {\hat X}\to X$ be the induced morphism. 
Hence by Theorem 7.6 and Definition--Proposition 1.5 
$$\split
e_{\orb}(x; X,D)&=\sum _{ z\in \{t\in {\hat X}\: \dim \varphi ^{-1}(t)>0\} }
e_{\orb} (z; {\hat X}, {\hat D})
-c_2({\hat \psi}, {\Olog {\hat \pi} {\hat X} {\hat D}})/\deg {\hat \pi}\\ 
&\le -{1\over 4}(c_1(\psi, K_{\hat X}+{\hat D}))^2.\\
\endsplit$$
But $c_1(\psi, K_{\hat X}+{\hat D})=-c_1(\psi , \psi ^{-1}B)$, so
$$e_{\orb}(x; X,D)\le -{1\over 4}(c_1(\psi, \psi ^{-1}B))^2.$$
Now the corollary follows from the following inequality
$$-(c_1(\psi, \psi ^{-1}B))^2\le -(\varphi ^*c_1(\psi, \psi ^{-1}B))^2-(c_1(\varphi, f ^{-1}B))^2=
-(c_1(f, f ^{-1}B))^2,$$
Q.E.D.

\proclaim{Corollary 7.8}
If $(\C ^2, D)$ is lc  then 
$$e_{\orb}(0; \C ^2,D)\le  \left( 1-{\mult _0 D \over 2}\right)^2.$$
In particular, $e_{\orb}(x; X,D)\le 1$ for any lc pair $(X,D)$.
\endproclaim

\pr
Let $f\: {\t X}\to X=\C^2$ be the minimal log resolution of $(\C^2 ,D)$ (or the blow up at $0$
if $(\C ^2, D)$ is smooth at $0$). Let $g \: {\hat X}\to X $ be the blow up at $0$ and
$h\: {\t X}\to {\hat X}$ the induced map.
Let $M_1$ and $M_2$ be two lines intersecting the $g$-exceptional divisor $F$ transversally
at single points, which are not blown up by $h$. Set $L_i=g(M_i)$ for $i=1,2$.

We can find rational numbers $0\le a_1\le a_2\le 1$ such that $(\C ^2, D+a_1L_1+a_2L_2)$
is lc but not lt. By construction $f^{-1}L_i=h^{-1}M_i=h^* M_i$ and $c_1 (f, h^* M_i)=-c_1 (f, h^*F)$.
Hence by Corollary 7.7
$$e_{\orb}(0; \C ^2,D)\le  -{1\over 4} ( c_1(f, -(a_1+a_2)h^*F))^2={1\over 4}(a_1+a_2)^2.$$
Since $\mult _0 (D+a_1L_1+a_2L_2 )\le 2$, we get the first inequality.
The second inequality follows from Lemma 7.1 and Corollary 7.5, Q.E.D.

\proclaim{Conjecture 7.9}
The function $f\: \C ^n \to \R$ sending $(a_1,\dots ,a_n)$ to $e_{\orb} (0; \C^2, \sum a_i D_i)$  
is continuous and nonnegative. Moreover, it is strictly decreasing for klt pairs.
\endproclaim

\proclaim{Problem 7.10}
Is the local orbifold Euler number semicontinuous in equisingular families?
Is it always semicontinuous?
\endproclaim

If the answers to these questions were positive then Corollary 7.8 would follow immediately
by deforming the divisor to the tangent cone 
and by Theorem 8.7 (cf. [Kol], the proof of Lemma 8.10).
In particular, this would imply Theorem 7.6. 
Indeed, if $e_{\orb} (0; \C ^2, L+D)> 0$ for some log canonical pair $(\C^2, L+D)$
then we can take a cyclic cover $\pi \: \C^2\to \C^2$ of degree $n$ branched along $L$. Let $M$ be the
ramification locus of $\pi$. Then $M$ is a line, $(\C^2, M+\pi ^*D)$ is lc and
$e_{\orb} (0; \C ^2, M+\pi^*D)=n \cdot e_{\orb} (0; \C ^2, L+D)$.
But $e_{orb}$ is bounded by Corollary 7.8, a contradiction.

Let $(X,D,x)$ be an analytic germ of a normal surface with $\Q$-boundary for which
one can define an orbifold Euler number. 
Let us recall that two pairs $(X_1,D_1,x_1)$ and $(X_2,D_2,x_2)$
are analytically  (topologically) equivalent if  and only if there exists
a biholomorphic map (a homeomorhism, respectively) $\varphi \: (X_1,x_1) \to (X_2,x_2)$ such that
$\varphi _*D_1=D_2$. 
By the definition $e_{\orb}(x;X,D)$ is an analytic
invariant, i.e., it depends only on the analytic type of $(X,D)$ at $x$.

According to Conjecture 7.9 we expect that
the log canonical threshold (see [Ko], 8.1) of a curve $C$ on a normal surface germ $(X,x)$
is given by the following formula
$$c(x, X,C)=\min \{ \alpha \:  e_{\orb}(x;X,\alpha C)=0 \} .$$
This invariant, in the surface case, can be read off from the embedded resolution
graph, so it is a topological invariant. 
However, Corollary 6.3 shows that an orbifold Euler number is not a topological
invariant of the pair $(X,\alpha C)$, since the Tjurina number changes under $\mu$-constant
deformations (see Theorem III.2.9.1 in Appendix by B. Teissier to [Z]). In fact, for smooth surface germs
$e_{\orb}(x;X,C)$ is a topological invariant of $(X,C)$ if and only if $C$ has a weighted
homogeneous singularity at $x$ (see Corollary 6.3).
It should be also true that if $C$ has a weighted homogeneous singularity at $x$ then
$e_{\orb}(x;X,\sum a_i C_i)$, where $C_i$ are irreducible components of $C$, is a topological
invariant of $(X,C)$ (cf. Theorems 8.3 and 8.7).

This shows that the computation of orbifold Euler numbers in general is a hopeless
problem and the most one can hope for is proving semicontinuity on
the moduli space of plane singularities with fixed semigroup and computing of $e_{\orb}$ at the generic points
of the main component of the moduli space. This problem is already difficult for a much simpler 
invariant, namely the Tjurina number (see [Z]) and the corresponding
algorithm is known only for irreducible curve singularities.

\medskip

\heading 8. Local orbifold Euler numbers for ordinary singularities
\endheading

In this section we compute local orbifold Euler numbers for ordinary singularities.
Since these numbers are analytic invariants it is sufficient to work with lines in $\C^2$.
First, we give a precise result for three lines in $\C ^2$ (see Theorem 8.3) and then a more general but 
slightly weaker result for any number of lines in $\C^2$ (see Theorem 8.7).
By Lemma 7.1 we can also deal with quotients of ordinary singularities (cf. Section 9 for 
quotients of ordinary singularities with three branches). 

\medskip

{\sl 8.1.} Let $L_1,\dots , L_n$ be $n$ distinct lines in $X=\C^2$ passing through $x=0$. 
Let $f\: ({\t X},E)\to (X,x)$ be the blow up of $X$ at  $x$.
Then $\psi \:{\t X}={\Bbb V}(\O _{\P ^1}(-1))\to \P^1$ is a geometric line bundle.
Assume $L_i$ is given by $f_i(x,y)=0$ in $X$ and set ${\t L}_i=f^{-1}L_i$ and $P_i=\psi ({\t L}_i)$.
$\t X$ is covered by two affine pieces $U_1$ and $U_2$:
$$f\: U_1=\C^2\to \C ^2, \quad f((z_1,t_1))=(z_1, z_1 t_1),$$
$$f\: U_2=\C^2\to \C ^2, \quad f((z_2,t_2))=(t_2z_2,  z_2).$$
The map $\psi$  is given by $\psi\: U_1\to \P ^1$,  $\psi((z_1,t_1))=[1, t_1],$
and $\psi\: U_2\to \P ^1$,  $\psi((z_2,t_2))=[t_2, 1].$

\proclaim{Lemma 8.2}
If $n\ge 1$ then 
$${\olog {\t X} {\sum _{i=1}^n{\t L}_i+E}}\simeq \psi ^* (\omega _{\P^1} (P_1+\dots +P_n)\op \O_{\P ^1}).$$ 
\endproclaim

\pr
We can assume that $f_1(x,y)=x$. Then ${\olog {\t X} {\sum {\t L}_i+E}}$ is generated by 
$$\omega_1={dz_1\over z_1}\quad \hbox{ and }\quad
\omega _2={dt_1\over {\prod _{i=1}^n f_i (1,t_1)}}$$
on $U_1$ and by
$$\omega_1'={dz_2\over z_2}+{dt_2\over t_2}\quad \hbox{ and }\quad
\omega _2'=-{dt_2\over {\prod _{i=1}^n f_i (t_2,1)}}$$
on $U_2$. Since 
$\omega _1|_{U_1\cap U_2}=\omega _1'|_{U_1\cap U_2}$ and 
$\omega _2|_{U_1\cap U_2}=t_2 ^{n-2}\omega _2'|_{U_1\cap U_2},$
we get the required isomorphism, Q.E.D.

\proclaim{Theorem 8.3}
Let $L_1$, $L_2$ and $L_3$ be 3 distinct lines in $\C ^2$ passing through  $0$. 
Set $D=a_1L_1+a_2L_2+a_3L_3$, where 
$0\le a_1\le a_2\le a_3\le 1$.
Then 
$$
e_{\orb} (0;\C^2,D)=\left\{ 
\aligned 0  &\quad\hbox{ if } a_1+a_2+a_3>2 \hbox{ (non lc case)},\\
          (1-a_1-a_2)(1-a_3)  &\quad \hbox{ if } a_3\ge a_1+a_2,\\
          {(a_1+a_2+a_3-2)^2\over 4}&\quad\hbox{ if } a_3<a_1+a_2\hbox{ and }
a_1+a_2+a_3\le 2.
\endaligned
\right.
$$
\endproclaim

\pr
We can choose coordinates in $(X,x)=(\C ^2,0)$ such that 
$L_1$, $L_2$ and $L_3$ are described by the equations 
$x=0$, $y=0$ and $x+y=0$, respectively. 
We can also find an integer $n$ such that
$a_i=1-l_i/n$ for some integers $l_i$. 
Let us consider a cone $Y$ in $\C^3$ over the curve $C\subset \P^2$ described 
by the equation $a^n+b^n+c^n=0$. 
Define the map $\pi\: Y\to X$ by  $\pi((a,b,c))=(a^n,b^n)$.  
This map is branched precisely over the three lines $L_i$ with
branching indices $n$. Let $M_i=(\pi^*L_i)_{red}$ be a reduced
scheme structure on the set theoretical inverse image of $L_i$. 
The divisor $M_i$ consists of $n$ lines intersecting in $y$.

The geometric line bundle 
$\varphi \:{\t Y}={\Bbb V}(\O _C(-1))\to C$ is a resolution of the cone 
singularity $Y$. Let $g\: {\t Y}\to Y$ be the contraction of 
the zero section $F$. Then the induced morphism 
${\t \pi} \:{\t Y}\to {\t X}$ is  finite.
Therefore to compute $e_{\orb} (x;X,D)$ it is sufficient to compute
$c_2(y,{\Olog {\t \pi} {\t X} {f^{-1}D+E} })$.

\proclaim{Lemma 8.4}
Set $\F_1=\O _C(-l_1)\op\O _C(-l_2)\op\O _C(-l_3)\op\O _C(n-l_1-l_2-l_3)$.
Let the map  $\F_1 \to \O _C\op\O _C(n)$ be given by the matrix
$$A={\pmatrix
a^{l_1}&b^{l_2} &c^{l_3} &0 \\
0 & -b^{l_2}c^n &  b^n c^{l_3}& a^{l_1} b^{l_2}c^{l_3} 
\endpmatrix}
$$
and let $\E$ be the image of this map.
Then $\Olog {\t \pi} {\t X} {f^{-1}D+E}=\varphi ^*\E$.
\endproclaim

\pr
The morphism $g\: {\t Y}\to Y$ is a blow up of the maximal ideal at $0$
and hence it is covered by three affine pieces $V_1$, $V_2$ and $V_3$.
Consider the morphisms
$$\split
&\sigma _1 \: V_1'=\C^3\to \C^3,\quad \sigma_1((u_1, v_1,w_1))=(u_1,u_1v_1,
u_1w_1),\\
&\sigma _2 \: V_2'=\C^3\to \C^3,\quad \sigma_2((u_2, v_2,w_2))=(u_2v_2,v_2,
v_2w_2),\\
&\sigma _3 \: V_3'=\C^3\to \C^3,\quad \sigma_3((u_3, v_3,w_3))=(u_3w_3,v_3w_3,
w_3).\\
\endsplit$$
Then $V_i=V_i'\cap {\t Y}$ and $g|_{V_i}={\sigma _i}|_{V_i}$.

Set  ${\t M}_i=g^{-1}M_i$.
By Lemma 8.2
$$\Omega _{\t Y}(\log {\t M}_1+{\t M}_2+{\t M}_3+F)={\t \pi}^*{\olog {\t X} {{\t L}_1+{\t L}_2+{\t L}_3+E}}
\simeq \O_{\t Y}\op \O_{\t Y}(-nF)$$
and we can use $\pi^* \omega_i$ and $\pi^*\omega _i'$ from the proof of Lemma 8.2 as local generators
of $\Omega _{\t Y}(\log {\t M}_1+{\t M}_2+{\t M}_3+F)$.
Using this basis one can see that
$\Olog {\t \pi} {\t X} {f^{-1}D+E}$ is given as a subsheaf of 
${\t \pi}^*{\olog {\t X} {{\t L}_1+{\t L}_2+{\t L}_3+E}}$
by the matrix $\varphi ^*A$. 
 
For example, on $V_1$ the sheaf  $\Olog {\t \pi} {\t X} {f^{-1}D+E}$ 
is generated by the forms 
$$\eta _1=
{d (v_1^n)\over v_1^{n-l_2}},\quad
\eta _2={d (w_1^n)\over w_1^{n-l_3}}\quad \hbox{and} \quad
\eta _3={d (u_1^n)\over u_1^n}.$$
The pull back of $\F_2 \to \O _C\op\O _C(n)$ to $\t Y$ is given on $V_1$ by the matrix
$$A_1={\pmatrix
1&v_1^{l_2} &w_1^{l_3} &0 \\
0 & -v_1^{l_2}w_1^n &  v_1^n w_1^{l_3}& v_1^{l_2}w_1^{l_3} 
\endpmatrix}
$$
so we need to check that the forms
$$\split
&{\t \pi }^*\omega _1=\eta_ 3, \quad
v_1^{l_2} {\t \pi }^*\omega _1-v_1^{l_2}w_1^n{\t \pi }^*\omega _2=\eta_1 +
v_1^{l_2}\eta_3,\\
&w_1^{l_3}{\t \pi }^*\omega _1-v_1^n w_1^{l_3}{\t \pi }^*\omega _2=
-\eta _2+w_1^{l_3}\eta _3\quad \hbox{and} \quad
v_1^{l_2}w_1^{l_3}{\t \pi }^*\omega _2=v_1^{l_2-n}\eta_2=w_1^{l_3-n}\eta_1 \\
\endsplit
$$
locally generate the sheaf ${\Olog {\t \pi} {\t X} {f^{-1}D+E}}$.
The last form belongs to this sheaf since on $V_1=(1+v_1^n+w_1^n=0)$ 
either $v_1\ne 0$ or $w_1\ne 0$ and the remaining forms generate  
$\Olog {\t \pi} {\t X} {f^{-1}D+E}$. 
Similarly one can check it on other affine pieces, Q.E.D.

\proclaim{Lemma 8.5}
Set $S=\C [a,b,c]$ and $R=S/(a^n+b^n+c^n)$ and assume that $0<l_i<n$ for $i=1,2,3$.
Then the complex
$$\dots \ra^{B'}R^4\ra^B R^4\ra^{B'}R^4\ra^B R^4\ra ^A\Im A\raa 0,$$
where
$$B={\pmatrix
0 &c^{l_3} &-b^{l_2} &-a^{n-l_1} \\
-c^{l_3} & 0 &  a^{l_1}& - b^{n-l_2} \\
b^{l_2} &  -a^{l_1}& 0& - c^{n-l_3} \\
a^{n-l_1} &  b^{n-l_2}& c^{n-l_3} & 0
\endpmatrix} ,
$$
$$B'={\pmatrix
0 &c^{n-l_3} &-b^{n-l_2} &-a^{l_1} \\
-c^{n-l_3} & 0 &  a^{n-l_1}& - b^{l_2} \\
b^{n-l_2} &  -a^{n-l_1}& 0& - c^{l_3} \\
a^{l_1} &  b^{l_2}& c^{l_3} & 0
\endpmatrix} ,
$$
is a minimal free resolution of $\Im A$ over $R$.
\endproclaim

\pr
Localizing at $(a,b,c)$ we can assume that $S$ and $R$ are local rings.
Note that $BB'=B'B=(a^n+b^n+c^n)\cdot \Id _{S^4}$ and all the matrix elements
of $B$ and $B'$ belong to the maximal ideal of $S$. Thus $(B,B')$
is a {\sl reduced matrix factorization} of $(a^n+b^n+c^n)$ over $S$ (see [E]).
In particular, by Corollary 6.3, [E], the complex
$$\dots \ra^{B'}R^4\ra^B R^4\ra^{B'}R^4\ra^B R^4\ra \coker B\raa 0$$
is a nontrivial periodic minimal free resolution over $R$ of  the maximal Cohen--Macaulay
$R$-module $\coker B$. 
So it is sufficient to prove that the natural map $\coker B \to \Im A$ is an isomorphism.

Since  $\det B=(a^n+b^n+c^n)^2$ and $(a^n+b^n+c^n)$ is a prime, $\coker B$ is a rank $2$
reflexive $R$-module by Proposition 5.6, [E].
Since $AB=0$ and both  $\coker B$ and $\Im A$ are rank $2$ reflexive modules, 
it is sufficient to check that $\coker B$ is a direct summand of $R^4$
at every height $1$ prime of $R$. This follows because at every height $1$ prime of $R$
the matrix factorization gives rise to the exact sequence
$$0\raa \coker B \ra R^4 \ra \coker B'\raa 0$$
splitting the surjection $R^4\to \coker B$, Q.E.D.

\medskip

Set $\F_2=\O _C(-l_2-l_3)\op \O _C(-l_1-l_3)\op\O_C (-l_1-l_2)\op\O _C(-n)$. 
As a corollary of  Lemma 8.5 we get that the complex of $\O _C$-modules
$$\F_2\ra ^B \F_1\ra \E \raa 0$$
is exact if $0<l_i<n$. In the remaining cases $\E$ is decomposable (and the complex
is still exact). It is easy to see that $e=\deg \E=n-l_1-l_2-l_3$.

Now Theorem 8.3 follows by a direct computation from  Theorem 1.10 and the following lemma:

\proclaim {Lemma 8.6}
Let $p=\max (-l_1,-l_2,-l_3,e)$. Then
$${\overline s} (\E)=\max \left(p,{ e\over 2}\right).$$
\endproclaim

\pr
It is clear from the definition of $\E$ that ${\overline s} (\E)\ge p$. 
Suppose that $e>2p$.
We need to prove that in this case $\E$ is semistable.

Let $\F_1'=\O _{\P ^2}(-l_1)\op\O _{\P ^2}(-l_2)\op\O _{\P ^2}(-l_3)\op\O 
_{\P ^2}(n-l_1-l_2-l_3)$ and
$\F_2'=\O _{\P ^2}(-l_2-l_3)\op \O _{\P ^2}(-l_1-l_3)\op\O_{\P ^2} 
(-l_1-l_2)\op\O _{\P ^2}(-n)$.

Consider the map of sheaves $\F_2'\to \F_1'$ given by the matrix $B$. 
Since $\det B=(a^n+b^n+c^n)^2$, this map is injective
and its cokernel $\G$ is supported on the curve $C$. 

Since the matrix $B$ is skew-symmetric one can easily see that
the kernel of the  map $\F_2\to \F_1$ is  isomorphic to $\E (-n)$
(here we also use equality $\E^* =\E (-\det \E)$ for a rank 2 vector bundle
$\E$). 

Hence we get the following commutative diagram:
$$\CD
@. @. 0@>>> 0 \\
@. @. @VVV @VVV \\
@. 0@>>> \F_2'(-C) @>>> \F_1' (-C) @>>> \G (-C) @>>> 0\\
@. @VVV @VVV @VVV @VVV @VVV \\
@. 0@>>> \F_2' @>>> \F_1'  @>>> \G  @>>> 0\\
@. @VVV @VVV @VVV @VVV @VVV \\
0@>>>\E(-n)@>>> \F_2 @>>> \F_1 @>>> \E @>>> 0\\
@. @. @VVV @VVV @VVV @. \\
@. @. 0 @>>> 0 @>>> 0 \\
\endCD
$$
Using the snake lemma we get an exact sequence
$$0\to \E(-n)\to \G(-n)\to  \G\to \E\to 0.$$
In particular, we get an inclusion $\E\hookrightarrow \G$.
Hence we have an inclusion $S^k\E\hookrightarrow S^k\G$ 
(using the fact that it is an inclusion at a general point of $C$
and $S^k\E$ is a torsion free $\O_C$-module). 
Note that $S^k \G$ is a subsheaf of $\H_k=\coker (S^k\F_2'\to S^k\F_1')$.

From the long exact cohomology sequence for the short exact sequence
$$0\to S^k\F_2'(m)\to S^k\F_1'(m)\to \H_k (m)\to 0$$
we see that
$$h^0(\P^2, \H_k (m))=0\quad \hbox{for}\quad m< -kp.$$
Since 
$$S^{2k}\E(-ke)\ot \O_C(-i)\hookrightarrow S^{2k}\G (-ke-i)
\hookrightarrow \H_{2k}(-ke-i)$$
and $ke\ge 2kp$ we get
$$h^0\left(
S^{2k}\E(-k\det \E)\ot \O_C(-i)\right)
=0 \quad\hbox{for}\quad i>0.$$
Hence by Corollary 1.9 the  vector bundle $\E$ is semistable, Q.E.D.

\proclaim{Theorem 8.7}
Let $L_1,\dots ,L_n$ be $n$ distinct lines in $\C ^2$ passing through  $0$. 
Set $D=\sum_{i=1}^na_iL_i$, where 
$0\le a_1\le a_2\le \dots \le a_n\le 1$, and $a=\sum_{i=1}^n a_i$.
Then 
$$
e_{\orb} (0;\C^2,D)=\left\{ 
\aligned 0  &\hbox{ if } a>2 \hbox{ (non lc case)},\\
          (1-a+a_n)(1-a_n)  &\hbox{ if } 2a_n\ge a,          
\endaligned
\right.
$$
and
$$e_{\orb} (0;\C^2,D)\le {\left(1-{a \over 2}\right)^2}\hbox{ if } 2a_n<a\le 2.$$
\endproclaim

\pr 
By increasing $n$ if necessary and adding to $D$ more lines with coefficient $0$,
we can assume that all $k_i = n a_i$ are integers. 

Let $\pi \: Y\to X=\C^2$ be a cyclic covering of order $n$ branched in $L_1+\dots +L_n$.
The blow up $g\:{\t Y}\to Y$ of $Y$ at $y=\pi ^{-1}(x)$ is the geometric line bundle $\varphi \:
{\t Y}={\Bbb V}(\O _C (-1))\to C$ over the curve $C$ given by  $z^n=f_1(x,y)\cdot \dots \cdot f_n(x,y)$
in $\P^2$. Hence we get a commutative diagram
$$\CD
C @<{\varphi}<< ({\t Y},F) @>{g}>> (Y,y)\\
@VV{p}V @VV{\t \pi}V  @VV{\pi }V \\
\P^1 @<{\psi}<< ({\t X},E) @>{f}>> (X,x).\\ 
\endCD
$$
Set $M_i= (\pi^* L_i)_{red}$, ${\t M}_i=g^{-1}M_i$, $Q_i =\varphi ({\t M}_i)$ and $l_i=n-k_i$.
To prove the theorem it is sufficient to compute ${c_2(g, \F)}=-n\cdot e_{\orb} (0;\C^2,D)$, where
$\F ={\Olog {\t \pi} {\t X} {\sum {k_i\over n}{\t L}_i+E} }$.
The computation follows easily from Theorem 1.10 and the following lemma.

\proclaim{Lemma 8.8}
Consider $\E =\F |_E$ as a bundle on $C=\varphi (E)$. 
Then $\F =\varphi ^* \E$ and $\E$ fits into the following
exact sequences:
$$0\raa \omega _C (-\sum (l_i-1) Q_i)\ra \E\ra \O _C\raa 0$$
and 
$$0\raa \O _C (-l_i Q_i)\ra \E\ra \omega _C (-\sum_{j\ne i} l_j Q_j)\raa 0$$
for $i=1,\dots , n$. Moreover,  the induced map
$$\omega _C (-\sum (l_i-1) Q_i)\op \bigoplus _{i=1}^n \O _C(-l_iQ_i)\to \E $$
is a surjection.
\endproclaim

\pr
By Lemma 8.2 
$${\olog  {\t Y} {\sum {\t M}_i+F} }={\t \pi}^*{\olog {\t X} {\sum _{i=1}^n{\t L}_i+E}}
\simeq \varphi ^* (\omega _{C} (Q_1+\dots +Q_n)\op \O_{C}),$$
so the lemma holds when all $l_i=0$. In general 
we just construct some maps leaving the details to the reader.

Let $a_i\in H^0(\O_C(Q_i))$ be the section corresponding to $Q_i$.
One can write it down explicitly as $z$ on $C-\bigcup _{j\ne i}Q_j$ and $1$ on $C-Q_i$.
Set $V_i=p^{-1}(\psi (U_i))$ for $i=1,2$. Then $(V_1, \prod a_i ^{l_i}p^*\omega _2)$
and $(V_2, \prod a_i ^{l_i}p^*\omega '_2)$ defines the map $\omega _C (-\sum (l_i-1) Q_i)\to \E$.

Let $\{ W_k\}$ be an affine cover of $C$ such that $\O_C(1)|_{W_k}$ is free, with generator $s_k$.
Then $s_k {\partial \over \partial s_k}=s_l{\partial \over \partial s_l}$  define a nowhere
vanishing section of ${{\widehat {{\t \pi} ^{*} \Der} _{\t X}(\log {\sum _{i=1}^n{k_i\over n}{\t L}_i+E})}}$
corresponding to the surjection $\E \to \O_C$.

The maps $\O _C (-l_i Q_i)\to \E$ are defined by $z^{-l_i}{p^* (df_i)}$ on $C-\bigcup _{j\ne i}Q_j$
and $p^* (df_i)$ on $C-Q_i$, Q.E.D.

\medskip

\re
Theorem 8.7 is a direct generalization of Theorem 8.3 except for the last case in which
we get inequality instead of equality. 
Equality, as in the proof of Theorem 8.3, is equivalent to
semistability of the vector bundle $\E$ from Lemma 8.8 (cf. Lemma 8.6 in the case of three lines).

\medskip

\heading  9. Quotients by unitary subgroups\endheading

The main aim of this section is a computation of local orbifold Euler numbers for log canonical
pairs with at most fractional boundary. Nontrivial examples of such pairs come from
quotients of $\C ^2$ by unitary subgroups of $GL (2, \C )$, whence the title of this section.

{\sl 9.1.} First, we need to introduce some notation. 
Let $n$ be a positive integer, $1\le q<n$ an integer coprime to $n$ and $\epsilon$
a primitive $n$th root of unity. Let us recall that the minimal resolution of the 
cyclic quotient singularity $(\a ^2, 0)/\Z _n$, $(x,y)\to (\epsilon x, \epsilon ^q y)$ 
consists of a chain of smooth rational curves $E_1, \dots ,E_s$ determined by the continued
fraction expansion
$${n\over q}=b_1-{1\over{b_2-{1\over {b_3-\dots }}}}.$$
We say that such chain is of type $\langle n, q\rangle$.
 $\langle 1, 0\rangle$ denotes an empty chain and $\langle 1,1 \rangle$ a single $(-1)$-curve. 

We say that the minimal log resolution of $(X,D,x)$ is of type  $\langle n, q; *_1\rangle$
($\langle n, q;*_1,*_2\rangle$) if its exceptional set is of type  $\langle n, q\rangle$
and $D$ has at most $1$ irreducible component (respectively: at most $2$ components)  
meeting the last (and the first) curve in the chain. 

The minimal log resolution of $(X,D,x)$ is of type  $\langle b;
\langle n, q; *_1\rangle,\langle n, q; *_2\rangle,\langle n, q; *_3\rangle \rangle$
if it is a star shaped tree of rational curves consisting of a central curve with
self intersection $-b$ and three chains $\langle n, q; *_1\rangle$ attached to it.

The type of the minimal log resolution depends only on the reduced pair
$(X,\lceil D\rceil)$ and in our notation  $*_i$ corresponds to the $i$th irreducible component of
$\Supp D$.

\medskip

{\sl 9.2.}
Let $(X,B=\sum b_i B_i, x)$ be a germ of a log canonical surface with a fractional
boundary, i.e., all the coefficients are of the form $b_i=1-{1\over m_i}$ for some 
$m_i \in \N \cup \infty$.  
In the above set up there is a natural interplay between unitary group actions on
the universal covering and the type of the resolution of singularity.
We will use it to compute all the local orbifold Euler
numbers for germs $(X,D,x)$, where $D=\sum d_i B_i$, $0\le d_i\le 1$.

Since we are interested mainly in the log canonical case we can assume that $(X, 0)$
is klt at $x$. Otherwise,  $e_{\orb} (x;X,D)=0$ by Corollary 7.5.

If $(X,0)$ is klt at $x$ then $(X, B, x)$ is analytically equivalent to
the quotient of $\C ^2$ by a finite subgroup $G\subset GL (2, \C)$ and the $B_i$
correspond to the components of the branch locus of the quotient map
$\pi\: \C ^2\to \C ^2/G$. The ramification index over a component with
coefficient $1-{1\over m}$ in $B$ is equal to $m$. 
Since $e_{\orb}$ depends only on the analytic type of the pair we can
assume that $X=\C ^2/G$.

Let ${\t G} \subset PGL (2,\C )$ be the projectivized group $G$, i.e., the
image of $G$ under the natural map $GL (2,\C )\to PGL (2,\C )$.

We consider the following two cases according to the type of 
the minimal log resolution of $(X,B)$:
\item{(1)} the minimal log resolution is of type $\langle n,q;*_1,*_2\rangle$,
\item{(2)} the minimal log resolution is of type
$\langle b; \langle n_1,q_1;*_1\rangle,\langle n_2,q_2;*_2\rangle, \langle n_3,q_3;*_3\rangle \rangle$.

In the first case $\t G$ is cyclic and in the second case $\t G$ is polyhedral
of type $\langle p_1,p_2,p_3\rangle$, where $p_i=n_im_i$. Let us recall that 
polyhedral groups are either dihedral $\langle 2,2,n\rangle$ or tetrahedral $\langle 2,3,3\rangle$ or octahedral 
$\langle 2,3,4 \rangle$ or icosahedral $\langle 2,3,5\rangle$.

\medskip
{\sl 9.3. Cyclic case}

\proclaim{Proposition 9.3.1}
Assume that the minimal log resolution of $(X,D=\sum d_iB_i,x)$ is of type $\langle n,q;*_1,*_2\rangle$.
Then
$$e_{\orb}(x;X,D)={(1-d_1)(1-d_2)\over n}.$$
\endproclaim

\pr
Since $(X,D,x)$ is a quotient of smooth pair $(\C^2, d_1(z_1=0)+d_2(z_2=0))$
by $\Z _n$ the proposition follows from Theorem 8.3 and Lemma 7.1, Q.E.D.

\medskip
\re
Proposition 9.3.1 was known; it is a local version of [Me2], Theorem 6.1; in the case $n=1$ it is also equivalent
to [Ti], Lemma 2.4.

\medskip
{\sl 9.4. Polyhedral case}

Let us recall Brieskorn's construction [Br].
Let $g\: ({\hat Y}, F)\to (Y=\C^2, 0)$ be the blow up of $\C^2$ at the origin.
The group $G$ acts on $\hat Y$. Let ${\hat \pi}\: {\hat Y}\to {\hat X} $ denote the quotient
map and let $E$ be the image of $F$. We have an induced proper morphism
$h\: ({\hat X}, E)\to (X,x)$, which is in fact obtained from the minimal resolution
$\t X \to X$  by contracting chains $\langle n_i, q_i\rangle$ to points.

\proclaim{Lemma 9.4.1} ([Me2], Theorem 2.9)
$\deg \pi =4s^2b_0$, where $b_0=b-\sum {q_i\over n_i}$ and $1+{1\over s}=\sum {1\over p_i}$.
\endproclaim

\pr
The action of $G$ on $\hat Y$ induces an action of $\t G$ on $F\simeq \P^1=\P (\C^2)$.
Hence the ramification index over $E$ is equal to ${|G|\over |{\t G}|}$ and 
${\hat \pi}^*E={|G|\over |{\t G}|}F.$
Since $c_2(g, {\hat \pi}^*(\O (E)\op \O (E)))= \deg {\hat \pi}\cdot c_2(h,\O (E)\op \O (E))$,
we have the equality
$$\left({|G|\over |{\t G}|}F\right)^2= \deg {\hat \pi}\cdot E^2.$$
But  $\deg {\hat \pi}=\deg \pi =|G|$, so
$\deg \pi =- |{\t G}|^2E^2$. It is easy to see that $-E^2=b_0$ and it is well known
that $|{\t G}|=2s$ (see [Yo], 11.2), which proves the lemma, Q.E.D.

\proclaim{Theorem 9.4.2}
Assume that the minimal log resolution of $(X,D=\sum d_iB_i,x)$ is of type 
$\langle b; \langle n_1,q_1;*_1\rangle,\langle n_2,q_2;*_2\rangle, \langle n_3,q_3;*_3\rangle \rangle$.
Set $\alpha = \sum {1-d_i\over n_i}$ and $\beta =\min \{ {1-d_i\over n_i}\}$.
Then
$$e_{\orb}(x;X,D)=\left\{ 
\aligned 0  &\quad\hbox{ if } \alpha < 1 \hbox{ (non lc case)},\\
          {{(\alpha -1)^2\over 4b_0}}&\quad\hbox{ if } 1\le \alpha < 2\beta +1,\\
          {(\alpha-1-\beta )\beta \over b_0}  &\quad\hbox{ if } \alpha\ge 2\beta +1.
\endaligned
\right. 
$$
\endproclaim

\pr
There exists a pair $(X,B=\sum b_iB_i,x)$ with  fractional boundary and such that the minimal
log resolution of $(X,B,x)$ is of the same type as that of $(X,D,x)$. We will use the notation
introduced in 9.2. 

Since $G$ is conjugate to a unitary subgroup, we can assume that $G\subset U(2)$. 
Let us recall that a unitary linear transformation is called {\sl a reflection}
if all but one  of its eigenvalues are equal to $1$. 
Let $N$ be a subgroup of $G$ generated by reflections.
Let $H$ be the maximal finite unitary reflection subgroup
of $GL (2, \C )$ containing $G$. 
Note that $|H|=4s^2$ by [Yo], 11.2. It is well known that $N\triangle G$ and $N\triangle H$. 
We have natural quotient maps
$\pi _1 \: Y=\C^2 \to Z=\C^2 /N$, $\pi _2\: Z \to X=Z/(G/N)$, $\varphi \: Y\to T=Y/H$
and $\psi \: Z\to T=Z/(H/N)$. The quotients $Z$ and $T$ are smooth by [Ch], Theorem A.
One can easily see that $\pi _2$ is ramified only at the point $x$ and 
$\pi _1$ is ramified only along irreducible curves $M_i={\pi_1}^{-1}B_i$
with ramification indices $m_i$.
The map $\varphi$ is ramified along lines $L_i=\psi (M_i)$ with ramification indices $p_i$.
Therefore $\psi$ is ramified along $L_i$ with ramification indices $n_i$ and
$\psi ^*L_i=n_i M_i$. 

It follows that
$${\Olog {\pi _2} X {D}}={\hlog Z {\sum d _iM_i}}={\Olog \psi T {\sum \left( 1-{1-d_i \over n_i}\right) L_i}}.$$
These equalities make sense even if the sheaves involved  are not well defined on $Z$; then 
we interpret them  on appropriate coverings. 
Therefore by Lemma 7.1 we get
$$e_{\orb } (x; X,D)\cdot \deg \pi_2= e_{\orb} \left(0;\C^2; \sum \left( 1-{1-d_i\over n_i}\right) L_i \right)
\cdot  
\deg \psi.$$
By Lemma 9.4.1 we have
$${\deg \psi \over \deg \pi_2}={|H/N|\over |G/N|}={|H|\over |G|}={4s^2\over 4s^2b_0}={1\over b_0}.$$
Hence the theorem follows by a simple computation from Theorem 8.1, Q.E.D.

\medskip
\re
The special case of Theorem 9.4.2 when $\langle p_1,p_2,p_3\rangle =\langle p_1,2,2\rangle $,
$d_2=b_2$ and $d_3=b_3$ is proved in [Me2], 6.1. In this case the computation can be reduced
 to a smooth log pair by covering techniques.

\medskip

\heading 10. Curves in surfaces of general type\endheading

The main result of this section is an effective version of Bogomolov's result
on boundedness of rational curves in surfaces of general type with
$c_1^2>c_2$ (see [Bo1]). Part (2) of Theorem 10.1 was known for curves containing only nodes:
see [Ti], Theorem 2.8 and [LM], Theorem 3. Its generalization to all curves was
conjectured by G. Tian in [Ti], 2.7--2.9.
Unfortunately, the author can prove this only for curves with ordinary singularities.
Theorem 10.1, (1) gives a result for all curves on surfaces of positive index.
It is also possible, by changing arguments slightly,
to get finiteness of rational curves
in the boundary cases (cf. [Ti], Theorem 2.9). 

Recently Miyaoka announced the proof of boundedness of canonical degree $K_SC$ in terms of
$c_1^2(S)$, $c_2(S)$ and $g(C)$  for any surface $S$ of general type (see [Mi3]).

\proclaim{Theorem 10.1}
Let $S$ be a surface of general type and let $C$ be a curve of geometric 
genus $g$ in $S$. 
\item{(1)} If $c_1^2>2c_2$  then
$$K_SC\le {(3c_2-c_1^2)(c_1^2+c_2)+\max (0,6(g-1)c_2) \over c_1^2-2c_2}.\leqno (10.1.1)$$
\item{(2)} If $c_1^2>c_2$ and $C$ has only ordinary singularities then
$$K_SC\le {3c_2-c_1^2+\max (0, 4g-4)\over c_1^2-c_2}c_1^2.\leqno (10.1.2)$$
\endproclaim

\pr
Let us fix $\alpha$ and let $f\: Y\to S$ be a morphism obtained by composing
blow ups of points at which the strict transform of $\alpha C$ is not log canonical 
and such that the strict transform of $C$ has maximal multiplicity among such points.
If $(S,\alpha C)$ is log canonical then $f$ is the identity map.

Let $n_P$ be the  number of blow ups of infinitely
near points lying over $P$ and let $E_P$ stand for $f^{-1}(P)$ with the reduced
structure (if $n_P\ne 0$). Let $F_{i,P}$ be the total transform on $Y$ of the divisor obtained by 
the $i$-th blow up 
of an infinitely near singular point $Q_{i,P}$ of $C$ lying over $P$.
Let $m_{i,P}$ be the multiplicity of $Q_{i,P}$ on the strict transform of $C$.
Then $K_Y=f^* K_S+\sum F_{i,P}$ and $f^*C=f^{-1}C+\sum m_{i,P}F_{i,P}$.

Let $E$ be the exceptional locus of $f$ with the reduced scheme structure. Then
$$\split
&c_2(Y, \alpha f^{-1}C)=e_{\top}(Y-(f^{-1}C\cup E))+(1-\alpha )e_{\top}(f^{-1}C-f^{-1}(\Sing C))+e_{\top}(E)\\
&+\sum _{y\in f^{-1}C\cap E}(e_{\orb}(y; Y, \alpha f^{-1}C)-e_{\orb}(y;Y, E))+
\sum _{P\in \Sing C -f(E)}e_{\orb} (P;S,\alpha C)\\
&=
\sum _{P\in f(E)}\left(e_{\top} (E_P) -(1-\alpha)+\sum _{y\in f^{-1}C\cap E_P}
(e_{\orb}(y; Y, \alpha f^{-1}C)-1)\right)\\
&+\sum _{P\in \Sing C-f(E)}\left(e_{\orb} (P;S,\alpha C)-(1-\alpha)\right)-\alpha e_{\top}(C)+c_2(S).\\
\endsplit
$$
Let $r_P$ be the number of analytic branches of $C$ passing through $P$. Using the equalities
$$e_{\top} (C)=2-2g-\sum _{P\in \Sing C}(r_P-1)$$ 
and
$$(K_Y+\alpha f^{-1}C)^2=(K_S+ \alpha C)^2+c_1^2(f, K_Y+\alpha f^{-1}C)$$
one can show that
$$\split
&3c_2(Y, \alpha f^{-1}C)-(K_Y+\alpha f^{-1}C)^2= 3c_2-c_1^2+3\alpha (2g-2)-2\alpha K_SC -\alpha ^2 C^2\\
&+\sum _{P\in \Sing C-f(E)}A_P+\sum _{P\in f(E)}( B_{P,1}+B_{P,2}),
\endsplit $$
where
$$A_P=3 (e_{\orb}(P;S,\alpha C)+r_P\alpha -1),$$
$$B_{P,1}=3\alpha r_P+3\sum _{y\in f^{-1}C\cap E_P}(e_{\orb}(y;Y, \alpha f^{-1}C)-1)$$
and
$$B_{P,2}=3e_{\top}(E_P)-3-\left(\sum_{i=1}^{n_P} (1-\alpha m_{i,P})F_{i,P}\right)^2.$$
Let $m_P$ denote the multiplicity of $C$ at $P$. Note that by Corollary 7.8 
$$A_P\le 3 \left( \left(1-{m_P\alpha \over 2}\right)^2+r_P\alpha -1\right)\le {3\over 4}
\alpha ^2m_P^2\le {3\over 2}\alpha ^2m_P(m_P-1).\leqno{(10.1.3)}$$
Let $m_y$ be the multiplicity of $f^{-1}C$ at $y\in  f^{-1}C\cap E_P$.
If $m_y=1$ then $e_{\orb}(y;Y, \alpha f^{-1}C)=1-\alpha$. If $m_y\ge 2$
then  $e_{\orb}(y;Y, \alpha f^{-1}C)\le (1-{\alpha m_y\over 2})^2$ by Corollary 7.8.
Hence 
$$B_{P,1}\le 3\alpha r_P-3\alpha \sum _{y\in f^{-1}C\cap E_P} m_y+{3\over 4}\alpha ^2
\sum _{m_y\ge 2}m_y ^2.$$
Since $\sum _{y\in f^{-1}C\cap E_P} m_y\ge r_P$ and 
$\sum _{m_y\ge 2}m_y ^2\le \sum _{m_y\ge 2}(m_y ^2+m_y(m_y-2))=\sum _{y\in f^{-1}C\cap E_P}2m_y(m_y-1)$
we get
$$B_{P,1}\le {3\over 2}\alpha ^2\sum _{y\in f^{-1}C\cap E_P}m_y(m_y-1).\leqno (10.1.4)$$
Since $F_{i,P}^2=-1$ and $F_{i,P}F_{j,P}=0$ for $i\ne j$ and $E_P$ is a tree of rational curves 
we get
$$B_{P,2}=3n_P+\sum _{i=1}^{n_p}(1-\alpha m_{i,P})^2=\sum _{i=1}^{n_p}(4-2\alpha m_{i,P}+
\alpha^2m_{i,P}^2).\leqno (10.1.5)$$
Since we blow up only non lc points we have $\alpha m_{i,P}>1$. Now we distinguish  two cases
depending on the assumptions.

\medskip

(1)  Fix a real number $\gamma$. If $\alpha \le {\gamma \over 3+\gamma}$ then
$m_{i,P}\ge m_{i,P}\alpha ({3\over \gamma}+1)>{3\over \gamma}+1$ and
$$4-2\alpha m_{i,P}+\alpha ^2m_{i,P}<3\alpha ^2m_{i,P}+ \alpha ^2(\gamma m_{i,P}-3-\gamma)m_{i,P}=
\alpha ^2(3+\gamma)m_{i,P}(m_{i,P}-1).$$
Therefore by (10.1.4) and (10.1.5) 
$$\split
B_{P,1}+B_{P,2}&\le {3\over 2}\alpha ^2\sum _{y\in f^{-1}C\cap E_P}m_y(m_y-1)+
 (3+\gamma)\alpha ^2\sum _i m_{i,P}(m_{i,P}-1)\\
&\le (3+\gamma)\alpha ^2\sum_{Q\to P}m_Q(m_Q-1),\\
\endsplit
$$
where the sum $\sum_{Q\to P}$ is taken over all infinitely near points of $C$ lying over $P$
(including $P$).

Since the pair $(Y, \alpha f^{-1}C)$ is log canonical Theorem 0.1 implies that
$$(K_Y+\alpha f^{-1}C)^2\le 3c_2(Y, \alpha f^{-1}C).$$
Therefore by (10.1.3) and by the above we have
$$\split
&2\alpha K_SC+\alpha ^2 C^2 - 3c_2+c_1^2 -3\alpha (2g-2)\le (3+\gamma)\alpha ^2
\sum _Q m_Q(m_Q-1)\\
&=(3+\gamma)\alpha ^2 \left((K_S+C)C-(2g-2)\right) \\
\endsplit
$$
for $\alpha \le {\gamma \over 3+\gamma}$, 
where the sum $\sum _Q$ is taken over all infinitely near points of $C$.

Set $x=K_SC$. Using the Hodge index theorem we have $C^2\le {x^2\over K_S^2}$. Substituting
into the above inequality we get
$$\beta \alpha ^2-2(x-3(g-1))\alpha+3c_2-c_1^2\ge 0,\leqno (10.1.6)$$
where $\beta= (2+\gamma)(x+{x^2\over K_S^2})+x-(2g-2)(3+\gamma )$.
If $\beta \le 0$ then substituting $\alpha ={\gamma \over 3+\gamma}$ into (10.1.6)
we get
$$x\le 3(g-1)+{3+\gamma\over 2\gamma } (3c_2-c_1^2).$$
Otherwise, set $\alpha _0={x-3(g-1)\over \beta}$. Then we can rewrite the above inequality
in the following form:
$$\beta (\alpha-\alpha_0)^2+3c_2-c_1^2\ge \beta \alpha_0^2$$
for $0\le \alpha \le {\gamma \over 3+\gamma}$.
We have three possibilities:
\item{(1.1)} either $\alpha _0\ge{\gamma \over 3+\gamma}$, or
\item{(1.2)} $\alpha _0<0$, in which case $x<3(g-1)$ or
\item{(1.3)} $0<\alpha _0\le{\gamma \over 3+\gamma}$. Then
we can set $\alpha =\alpha _0$ and then we get $\beta \alpha_0^2\le 3c_2-c_1^2$.

If we set  $\gamma={c_1^2\over c_2}-2$  then in all the cases we get the inequality (10.1.1)
by straightforward computations.

\medskip

(2) By construction $m_{Q_{n_P}}\le m_{Q_{n_P-1}}\le\dots \le m_{Q_1}=m_P$. So if the strict
transform of $C$ has an ordinary singularity at $Q_{n_P}$ then $\alpha m_{i,P}\ge \alpha m_{n_P,P}>2$
and $m_{i,P}\ge 3$. 
(Otherwise, the strict transform of $\alpha C$ is lc at $Q_{n_P}$, a contradiction.)
Hence  in this case
$$B_{P,2}\le \alpha^2\sum _{i=1}^{n_p}m_{i,P}^2+ {1\over 2}\alpha ^2 \sum m_{i,P}(m_{i,P}-3)=
{3\over 2}\alpha ^2\sum m_{i,P}(m_{i,P}-1).$$
As in (1) this, together with (10.1.3), yields
$$2\alpha K_SC+\alpha ^2 C^2 - 3c_2+c_1^2 -3\alpha (2g-2)\le
{3\over 2}\alpha ^2 \left((K_S+C)C-(2g-2)\right) .$$
As in (1) we can use the Hodge index theorem.
Then the inequality (10.1.2) follows by straightforward computations, Q.E.D.

\medskip
\re It should be possible to improve the above method to bound $K_SC$ for all curves on surfaces with 
$c_1^2>c_2$. It is easy to give such a bound in terms of $c_1^2, c_2, g$ and the maximal multiplicity
of points occuring on $C$.

\medskip

\heading 11. Plane curves \endheading

{\sl 11.1.} Let $C$ be a reduced curve on a smooth surface $X$. Then
$$e_{\top}(C)=-(K_X+C)C+\sum _{P\in \Sing C}\mu _P,$$
where $\mu _P$ is a Milnor number of $C$ at $P$.
Therefore, if $(X, \alpha C)$ is log canonical and $K_X+\alpha C$
is pseudoeffective then we can rewrite the inequality in Corollary 5.2
in the following form:
$$\sum _{P\in \Sing C}3(\alpha(\mu _P-1)+1-e_{\orb}(P;X, \alpha C))\le
3c_2-c_1^2 +\alpha K_XC+(3\alpha -\alpha ^2)C^2. \leqno (11.1.1)$$

\medskip

{\sl 11.2. Curves with a maximal number of cusps.}

Theorem 9.4.2 can be applied to study a large class of singularities, 
including all simple curve singularities. For example, if $C$ has an ordinary
cusp $x^2=y^3$ at $0\in \C^2$ then
$$e_{\orb}(0;\C^2,\alpha C)=\left\{ 
\aligned  {1-2\alpha}  &\hbox{ if } 0\le \alpha\le {1\over 6},\\          
          {3\over 2}{{\left(\alpha -{5\over 6}\right)^2}}&\hbox{ if } {1\over 6}\le \alpha \le {5\over 6},\\
          0  &\hbox{ if } {5\over 6}\le \alpha \le 1\hbox{ (non lc case)}.
\endaligned
\right. 
$$
If we apply this result to the pair $(\P^2, {\sqrt 73 -1\over 24}C)$ then we get the following
bound on the maximal number $s(d)$ of ordinary cusps on a plane curve 
(with simple singularities only) of degree $d$:
$$\limsup _{d\to \infty}{s(d)\over d^2}\le {125+\sqrt 73 \over 432}.\leqno (11.2.1)$$ 
This problem was considered by Hirzebruch in [Hr2], Section 8, where one can find the bound ${5\over 16}$. 
Hirano [Hn], Corollary 3, exhibited examples with
$$\limsup _{d\to \infty}{s(d)\over d^2}\ge {9\over 32},$$
but the author's hope is that (11.2.1) is optimal.

\medskip

{\sl 11.3. Arrangements of lines.}

An {\sl arrangement of $k$ lines} is a set of $k$ distinct lines in $\P ^2$.
Let $t_r$ be the number of points lying on exactly $r$ lines of the
arrangement.

\proclaim{Proposition 11.3.1}
If $t_r=0$ for $r>{2\over 3}k$ (i.e., the arrangement ``does not contain large pencils'') then
$$\sum_{r\ge 2}rt_r\ge\left\lceil{1\over 3}k^2+k\right\rceil\quad\hbox{and}\quad
\sum_{r\ge 2} r^2t_r\ge \left\lceil{4\over 3}k^2\right\rceil.\leqno (11.3.1.1)$$
\endproclaim

\pr
Let $C$ be a sum of lines in the arrangement.
By assumption there exists $\alpha$ such that $K_{\P^2}+\alpha C$ is nef and $(\P^2,\alpha C)$ is log 
canonical. If we apply inequality (11.1.1) and Theorem 8.7 to any such $\alpha$ and use
$\sum_{r\ge 2}t_rr(r-1)=k(k-1)$ then we get required inequalities, Q.E.D.

\medskip
{\sl 11.3.2. Example.}

Equality  holds in (11.3.1.1) for almost all arrangements coming from reflection groups
(see [Hr1], (1.2)). The only arrangements for which equality fails are some real arrangements.

For example there exists an  arrangement of lines $A^0_{3m}$ coming from a unitary reflection group  for which
$k=3m$, $t_3=m^2$, $t_m=3$ and $t_r=0$ for all other $r$. Hence we have infinitely many examples
for which equality holds in (11.3.1.1).

\enddocument